\documentclass[a4paper,11pt,freqn,leqno]{article}
\usepackage{amssymb,amsmath}
\def\Line#1{\hbox to\textwidth{#1}}%

\newtheorem{theorem}{Theorem}[section]
\newtheorem{lemma}{Lemma}[section]
\newtheorem{definition}{Definition}[section]
\newtheorem{cor}{Corollary}[section]
\newtheorem{prop}{Proposition}[section]

\newcommand\qed{{\unskip\nobreak\hfil\penalty50\hskip2em\vadjust{}
    \nobreak\hfil$\Box$\parfillskip=0pt\finalhyphendemerits=0\par}}

\renewcommand{\a}{\alpha}
\renewcommand{\b}{\beta}
\renewcommand{\d}{\delta}
\newcommand{\D}{\Delta}
\newcommand{\e}{\varepsilon}
\newcommand{\g}{\gamma}
\newcommand{\G}{\Gamma}
\renewcommand{\l}{\lambda}
\renewcommand{\L}{\Lambda}

\renewcommand{\O}{\Omega}
\renewcommand{\o}{\omega}

\renewcommand{\>}{\rangle}

\renewcommand{\i}{\infty}

\newcommand{\bE}{{\mathbb E}}
\newcommand{\bN}{{\mathbb N}}

\newcommand{\bR}{{\mathbb R}}
\newcommand{\bL}{{\mathbb L}}
\newcommand{\R}{{\mathbb R}}

\newcommand{\cA}{{\mathcal A}}
\newcommand{\cM}{{\mathcal M}}
\newcommand{\cN}{{\mathcal N}}
\newcommand{\cP}{{\mathcal P}}
\newcommand{\cF}{{\mathcal F}}
\newcommand{\cH}{{\mathcal H}}

\newcommand{\cX}{{\mathcal X}}
\newcommand{\cI}{{\mathcal I}}

\begin{document}

\title{Coagulation, Diffusion and the Continuous Smoluchowski Equation}
\author{Mohammad Reza Yaghouti \\
Amirkabir University\\ Mathematics and Computer Science Faculty \\
and \\
Fraydoun Rezakhanlou\thanks{This work is supported in part by NSF
grant DMS-0707890.}
\\
Mathematics Department, UC Berkeley \\
and \\
Alan Hammond \\
Courant Institute,  New York University}

\maketitle \setcounter{equation}{0}
\noindent{\bf Abstract.} 
The Smoluchowski equation is a system of partial differential equations modelling the diffusion and binary coagulation of a large 
collection of tiny particles. 
The mass parameter may be indexed either 
by positive integers, or by positive reals, these corresponding to the discrete or the continuous 
form of the 
equations.
In dimension $d \geq 3$, 
we derive the continuous Smoluchowski PDE as a kinetic limit of a microscopic model of Brownian particles liable to coalesce, using a similar 
method to that used to derive the discrete form of the equations in \cite{HR1}. The principal innovation is a correlation-type bound on 
particle locations that permits the derivation in the continuous context while simplifying the arguments of \cite{HR1}.
We also comment on the scaling satisfied by the continuous Smoluchowski PDE, and its potential implications for blow-up of solutions of 
the equations.
 
\begin{section}{Introduction}
 It is a common practice in statistical mechanics to
formulate a microscopic model with simple dynamical rules in order
to study a phenomenon of interest. 
In a colloid, a population of comparatively massive particles is agitated by the bombardment of much smaller particles in the ambient 
environment: the motion of the colloidal particles may then be modelled by Brownian motion.
Smoluchowski's equation provides a macroscopic description for the
evolution of the cluster densities in a colloid whose particles are prone to binary coagulation. Smoluchowski's
equation comes in two flavours: discrete and continuous. In the
discrete version, the cluster mass may take values in the set of positive integers,
whereas, in the continuous version, the cluster mass take values in
$\mathbb{R}^{+}$. 
Writing $f_{n}(x,t)$ for the density of clusters
(or particles) of size $n$, this density evolves according to
\begin{equation}\label{eq1.1}
\frac{\partial{f}_{n}}{\partial{t}}=d(n)\triangle{f}_{n}(x,t)
+Q_{+}^{n}(f)(x,t)-Q_{-}^{n}(f)(x,t),
\end{equation}
where
\begin{equation}\label{eq1.2}
Q_{+}^{n}(f)=\int_{0}^{n}\beta(m,n-m)f_{m}f_{n-m} dm,
\end{equation}
and
\begin{equation}\label{eq1.3}
Q_{-}^{n}(f)=2\int_{0}^{\infty}\beta(m,n)f_{m}f_{n} dm,
\end{equation}
in the case of the continuous Smoluchowski equation. In the
discrete case, the integrations in ~\eqref{eq1.2} and \eqref{eq1.3}
are replaced with summations.

In~\cite{HR1} and \cite{HR2}, we derived the discrete Smoluchowski
equation as a many particle limit of a microscopic model of
coagulating Brownian particles. (See also \cite{Lang}, \cite{N} and \cite{GKO} for similar results.) The main
purpose of the present article is the derivation of (1.1) in the
continuous case. 
We introduce a simpler approach to that used in 
\cite{HR1} and \cite{HR2}.
We will present a robust argument that allows us to
circumvent some induction-based steps of~\cite{HR1} and
\cite{HR2} (which anyway could not be applied in the continuous case).
As such, an auxiliary purpose of this article is to present a shorter proof of the kinetic limit derivations of Smoluchowski's equation
 given in  
\cite{HR1} and \cite{HR2}.
The
main technical tool is a correlation-type bound on the particle
distribution that seems to be applicable to general systems of
Brownian particles. To explain this further, we need to sketch the
derivation of Smoluchowski's equation and explain the essential role
of the correlation bounds.

The microscopic model we study in this article consists of a large
number of particles which move according to independent Brownian
motions whose diffusion rates $2d(m)$ depend on their mass $m \in
(0,\infty)$. Any pair of particles that approach to within a certain
range of interaction are liable to coagulate, at which time, they
disappear from the system, to be replaced by a particle whose mass
is equal to the sum of the masses of the colliding particles, and whose
location is a specific point in the vicinity of the location of the coagulation. 
This range of interaction is taken to be equal to a
parameter $\epsilon$, whose dependence on the mean initial total number $N$
of particles is given by $N=k_{\epsilon}Z$ for a constant $Z$,
where
\begin{displaymath} k_{\epsilon} = \left\{ \begin{array}{ll}
\epsilon^{2-d} & \textrm{\ \ \ if $d\geq 3,$}\\
|\log\epsilon| & \textrm{\ \ \ if $d=2$}.\\
\end{array} \right.
\end{displaymath}
This choice will ensure that a particle experiences an expected
number of coagulations in a given unit of time that remains bounded
away from zero and infinity as $N$ is taken to be high.

Our main result is conveniently expressed in terms of empirical
measures on the locations $x_{i}(t)$ and the masses $m_{i}(t)$ of
particles. We write $g(dx,dn,t)$ for the measure on
$\mathbb{R}^{d}\times[0,\infty)$ given by
\begin{equation} \nonumber
g(dx,dn,t)=k_{\epsilon}^{-1}
\sum_{i}\delta_{(x_{i}(t),m_{i}(t))}(dx,dn).
\end{equation}
Our goal is to show that, in the low $\epsilon$ limit, the measure $g$
converges to $f_{n}(x,t)dx~dn$, where $f_n$ solves the system
(\ref{eq1.1}). 
The main step in the proof requires
the replacement of the microscopic coagulation propensity $\a(n,m)$ (that we will shortly describe precisely)
of particles of
masses $n$ and $m$ with its macroscopic analogue $\b(n,m)$. The main technical
tool for this is a correlation bound which reads as follows, in the case that the coefficient $d(m)$ is non-increasing in $m$:
\begin{eqnarray} \label{eq1.4}
&& \mathbb{E}\int_0^\i\sum_{i_{1},...,i_{k}}
K\big(x_{i_{1}}(t),..,x_{i_{k}}(t)\big) \prod^{k}_{r=1}
d\big(m_{i_{r}}(t)\big)^{\frac{d}{2}}m_{i_{r}}(t) dt \\ \nonumber & &
\ \ \ \ ~~~\leq const.~~\mathbb{E} \sum_{i_{1},...,i_{k}}
\hat K\big(x_{i_{1}}(0),m_{i_1}(0),...,x_{i_{k}}(0),m_{i_1}(0)\big)
\prod^{k}_{r=1}d\big(m_{i_{r}}(0)\big)^{\frac{d}{2}} m_{i_{r}}(0).
\end{eqnarray}
Here, $\mathbb{E}$ denotes the expectation with respect to the
underlying randomness,   $K:(\mathbb{R}^{d})^{k}\rightarrow
\mathbb{R}$ is any non-negative bounded continuous function,
and $\hat K=-\Big(d(m_{i_1})\triangle_{x_{i_1}}+\dots+d(m_{i_k})\triangle_{x_{i_k}}\Big)^{-1}~K $. We refer to Section 4 for 
the corresponding correlation 
inequality when
the function $d(\cdot)$ is not non-increasing.

In fact, we need \eqref{eq1.4} only for certain examples of $K$ with
$k=2,3$ and $4$. It was these examples that were treated
in~\cite{HR1} and \cite{HR2} with rather ad-hoc arguments based on
an inductive procedure on the mass of the particles. Those arguments
seem to be specific to the discrete case and cannot be generalized
to the continuous setting.  Moreover, the bound \eqref{eq1.4}
implies that the macroscopic particle densities belong to $L^{p}$ for
given $p\ge 2$, provided that a similar bound is valid initially. This rather straightforward consequence of
\eqref{eq1.4} is crucial for the derivation of the macroscopic
equation. The corresponding step in~\cite{HR1} and \cite{HR2} is
also carried out with a method that is very specific to the discrete
case and does not apply to the continuous setting. This
important consequence of \eqref{eq1.4} simplifies the proof
drastically and renders the whole of section 4 of~\cite{HR1} redundant.

We state and prove our results when the dimension is at least three.
However, our proof for the correlation bound \eqref{eq1.4} works in any
dimension, and an interested reader may readily check that, as in this
article, the approach of~\cite{HR2} may be modified to establish
Theorem~1.1 in dimension two.

We continue with the description of the microscopic model and the
statement of the main result.

As a matter of convenience, we introduce two different microscopic models, that differ only in whether 
the number of particles is initially deterministic or random. We will refer to the model as deterministic
 or random accordingly. In either case, 
we define a sequence of microscopic models, indexed 
by a postive integer $N$.

A countable set $I$ of symbols is provided. A configuration
$\mathbf{q}$ is an $\mathbb{R}^{d} \times (0,\infty)$-valued
function on a finite subset $I_{\mathbf{q}}$ of $I$. For any $i \in
I_{\mathbf{q}}$, the component $ q(i)$ may be written as
$(x_{i},m_{i}).$ The particle labelled by $i$ has mass $m_{i}$ and
location $x_{i}$.  

In the deterministic case, the index $N$ of the model specifies the total number of particles present at time zero.
 Their placement is given as follows.
 There is a given function $h:\mathbb{R}^{d}\times (0,\infty)  \rightarrow
[0,\infty)$, with $h_n(x):=h(x,n)$, where $\int_0^\i \int_{\R^d} h(x,n)dxdn<\i$.
We set $Z=\int^{\infty}_{0} \int_{\mathbb{R}^{d}} h_{n}(x)\ dx
dn \in (0,\infty)$ and choose $N$ points in $(0,\infty) \times
\mathbb{R}^{d}$ independently according to a law whose density at
$(x,n)$ is equal to $h_{n}(x)/Z$. Selecting arbitrarily a set of $N$
symbols $\{i_{j}:j \in \{1, \ldots,N\}\}$ from $I$, we define the
initial configuration $\mathbf q (0)$ by insisting that $q_{i_{j}}(0)$ is
equal to the $j$-th of the randomly chosen members of $(0,\infty)
\times \mathbb{R}^{d}.$

In the random case, the index $N$ gives the mean number of initial particles. We suppose given some measure $\gamma_N$
 on positive integers that satisfies $\mathbb{E} \big( \gamma_N \big) = N$ and ${\rm Var} \big( \gamma_N \big) = o(N^2)$. 
The initial particle number, written $\cN$, is a sample of $\gamma_N$. The particles present at time zero are scattered in 
the same way as they are in the deterministic case. The subsequent evolution, whose randomness is independent of the sampling of $\cN$, 
is also the same as in the deterministic setting.

To describe this dynamics, set a parameter $\epsilon > 0$
according to $N = k_{\epsilon} Z$, as earlier described. 
Let $F:\{\mathbb{R}^{d} \times
[0,\infty)\}^{I} \rightarrow [0,\infty)$ denote a smooth function,
where its domain is given the product topology. The action on $F$ of
the infinitesimal generator $\mathbb{L}$ is given by
\begin{equation}\nonumber
(\mathbb{L}F)(\mathbf{q})=\mathbb{A}_{0}F(\mathbf{q})+\mathbb{A}_{c}F(\mathbf
{q}),
\end{equation}
where the diffusion and collision operators are given by
\begin{equation} \nonumber
\mathbb{A}_{0}F(\mathbf{q})=\sum_{i \in I_{\mathbf{q}} }
d(m_{i})\triangle_{x_{i}}F
\end{equation}
and
\begin{eqnarray} \label{eq1.5}
 \mathbb{A}_{c}F(q)&=&\sum_{i,j \in I_{q}} \epsilon^{-2}
V\Big(\frac{x_{i}-x_{j}}{\epsilon} \Big) \alpha(m_{i},m_{j}) \\
\nonumber  &&\Big[\frac{m_{i}}{m_{i}+m_{j}}
F(S^{1}_{i,j}\mathbf{q})+\frac{m_{j}}{m_{i}+m_{j}}
F(S^{2}_{i,j}\mathbf{q})-F(\mathbf{q}) \Big].
\end{eqnarray}
Note that:
\begin{itemize}
\item  the function $V:\mathbb{R}^{d} \rightarrow [0,\infty)$ is
assumed to be H\"older continuous, of compact support, and with
$\int_{\mathbb{R}^{d}} V(x)dx=1$.
\item we denote by $S^1_{i,j}\mathbf{q}$ that configuration formed from 
$\mathbf{q}$ by
removing the indices $i$ and $j$ from $I_{\mathbf{q}}$, and adding a
new index from $I$ to which  $S^1_{i,j}\mathbf{q}$ assigns the value
$(x_i,m_i + m_j)$. The configuration  $S^2_{i,j}\mathbf{q}$ is
defined in the same way, except that it assigns the value  $(x_j,m_i
+ m_j)$ to the new index. The specifics of the collision event then
are that the new particle appears in one of the locations of the two
particles being removed, with the choice being made randomly with
weights proportional to the mass of the two colliding particles.
\end{itemize}

\noindent
\textbf {Convention.} Unless stated otherwise, we will adopt a notation whereby all the
index labels appearing in sums should be taken to be distinct.

 We refer the reader to \cite{HR1} and \cite{R} for the reasons for choosing
 $N=\epsilon^{d-2}Z$, the form of
 the collision
term in (1.5), and the interpretations of the various terms.

Let us write $\mathcal{M}_Z(\mathbb{R}^{d}\times[0,\infty))$ for the space of
non-negative measures $\pi$ on $ \mathbb{R}^{d}\times[0,\infty)$ such that 
\[
\pi\left(\mathbb{R}^{d}\times[0,\infty)\right)\le Z.
\]
This space is equipped with the topology of vague convergence which turns $\cM_Z$ into a compact metric space. 
We also write $\mathcal{M}_Z \Big(  \mathbb{R}^{d}\times[0,\infty)^2\Big) $ for the space of non-negative measures $\mu$ such that
for every positive $T$, $\mu\Big(  \mathbb{R}^{d}\times[0,\infty)\times[0,T]\Big) \le TZ$,
 which is also compact with respect to the topology of vague convergence. 
This space has a closed subspace $\cX$ which consists of measures $\mu$ such that
$\mu\Big(  \mathbb{R}^{d}\times[0,\infty)\times[t_1,t_2]\Big) \le (t_2-t_1)Z$, for every $t_1\le t_2$. 
As we will show in Lemma~6.2 of Section~6, the space $\cX$ consists of
 measures $\mu(dx,dn,dt)=g(dx,dn,t)dt$ with $t\mapsto g(dx,dn,t)$ a Borel-measurable function from
$[0,\i)$ to $\mathcal{M}_Z(\mathbb{R}^{d}\times[0,\infty))$.
We will denote by $\mathbb{P}_{N}=\mathbb{P}^{\epsilon}$ the probability measure
on functions from $t \in [0,\infty)$ to the configurations
determined by the process at time $t$. Its expectation will be
denoted $\mathbb{E}_{N}$. Setting
\begin{equation} \nonumber
g ^{\epsilon} (dx,dn,t) = \epsilon^{2-d} \sum_{i}
\delta_{(x_{i}(t),m_{i}(t))}(dx,dn),
\end{equation}
the law of
\begin{equation}\nonumber
\mathbf{q} \mapsto g^{\epsilon} (dx,dn,t)dt
\end{equation}
with respect to $\mathbb{P}^{\epsilon}$ induces a probability
measure $\mathcal{P}^{\epsilon}$ on the space $\cX$. We note that, since the space $\cX$ is a compact metric space, 
the sequence $\cP_\e$ is precompact with respect to the topology of weak convergence.

For the main result of this article, we need the following assumptions on $\a(\cdot,\cdot)$ and $d(\cdot)$:

\noindent
\textbf{Hypothesis 1.1.}
\begin{itemize}
\item The diffusion coefficient $d:(0,\infty)\to (0,\infty)$
is a bounded continuous function and there exists a uniformly positive continuous function $\phi:(0,\i)\to(0,\i)$ such that 
both $\phi(\cdot)$ and $\phi(\cdot)d(\cdot)$ are non-increasing.
\item The function $\a:(0,\infty)\times (0,\i)\to (0,\infty)$ is a bounded symmetric continuous function satisfying
\[
\sup_{n\le L}\sup_m\frac {\a(n,m)}{md(m)^{\frac d2}\phi(m)^{d-1}}<\i,
\]
for every $L>0$.
\end{itemize}

\noindent
{\textbf {Remarks.}}
\begin{itemize}
\item
 The condition that the function $\phi:(0,\infty) \to (0,\infty)$ exist is rather mild and is satisfied if $d(\cdot)$
is non-increasing. This condition requires that heavier particles
to diffuse slower which is natural from a physical point of view.
In fact when  $d(\cdot)$ is non-increasing, then we can simply choose $\phi(m)\equiv 1$. 
Also, if $d(\cdot)$ is non-decreasing, then the function $\phi$ exists and  can be chosen to be $\phi(m)=d(m)^{-1}$.
 From these two cases, we guess that 
the first condition is related to the variation of the function $d(\cdot)$.
As we will show in Lemma 2.2 of Section 2, the existence of such a function $\phi$ is equivalent to assuming that the total 
negative variation of 
$\log d(\cdot)$ 
over each interval $[n,\i)$, $n>0$, is finite.
\item
We note that if  the function $d(\cdot)$ is non-increasing,
 then the second condition for small $m$ and $n$ is equivalent to saying that  $\a(m,n)\le C\min(m,n)$.
However, when $m$ and $n$ are large,  the second condition is satisfied if for example
$\a(m,n)\le C md(m)^{d/2}nd(n)^{d/2}$. In summary, the second condition is rather mild if $m$ and $n$ are large,
but much more restrictive if both are
small.   
 Our stipulation that
$d$ be bounded is more restrictive in the case for values of its argument
close to zero, since it is reasonable to assume that very
light particles diffuse rapidly.
 \end{itemize}

We also need the following assumptions on the initial data $h$:

\noindent
\textbf {Hypothesis 1.2.}
\begin{itemize}
\item
$\int_0^\i  \int h_n(x)dxdn<\i$.
\item
 ${\bar h}_k*\l_k\in
L^{\infty}_{loc}(\mathbb{R}^{d})$, for $k=2,3$ and $4$,
 where $\bar h_k=\int_0^\i nd(n)^{\frac d2-\frac 1k}\phi(n)^{\frac {dk}2-1}h_n\ dn$ and
$\l_k(x) = |x|^{\frac 2k-d}$.
\item
\[
\int\hat h(x)\hat h(y)|x-y|^{2-d}dxdy<\i
\]
where $\hat h=\int_0^\i (n+1)h_ndn$.
\end{itemize}

\noindent
{\textbf {Remark.}} Recall that if $d(\cdot)$ is non-increasing, then we may choose $\phi=1$. In this case, Hypothesis 1.2
is satisfied if $\hat h\in L^1\cap L^\i$.

To prepare for the statement of our main result, 
we now recall the weak formulation of the system (\ref{eq1.1}).
Firstly, recall that a non-negative measurable function
$f: \R^d \times [0,\infty) \times [0,\infty) \to [0,\infty)$
is a {\textbf {weak solution}} of (1.1) subject to the initial condition $f(x,n,0)=
h_n(x)$, if for every smooth function $J:\bR^d\times (0,\i)\times [0,\i)\to\bR$ of compact support,
\begin{eqnarray*}
\int_0^\i \int_{\R^d} f(x,n,t)J(x,n,t)dxdndt&=&\int_0^\i\int h_n(x)J(x,n,0)dxdn\\
&&+\int_0^t\int_0^\i \int_{\R^d}\frac{\partial J}{\partial t}(x,n,s)f(x,n,s)dxdnds\\
&&+\int_0^t\int_0^\i \int_{\R^d} d(n)\triangle J(x,n,s)f(x,n,s)dxdnds\\
&&+\int_0^t\int_0^\i \int_0^\i \int_{\R^d} \b(m,n)f(x,n,s)f(x,m,s)\\
&&\quad\quad\quad \tilde J(x,m,n,s)dxdndmds,
\end{eqnarray*}
where
\[
\tilde J(x,m,n,s)=J(x,m+n,s)-J(x,m,s)-J(x,n,s).
\]
Following Norris [9], we define an analagous measure-valued notion of weak solution.
\begin{definition}
Let us write $M[0,\i)$ for the space of non-negative measures on the interval $[0,\i)$. We equip this space with the topology of 
vague convergenece. A measurable function $f: \R^d \times [0,\infty)\to M[0,\i)$ 
is called a measure-valued weak solution of (\ref{eq1.1})
if, firstly, for each $\ell > 0$, 
the functions $g_\ell, h_\ell\in L^1_{loc}$, where
\[
g_\ell(x,t)=\int_0^\ell f(x,t,dn),\ \ \ h_\ell(x,t)=\int_0^\i\int_0^\ell\b(m,n)f(x,t,dn)f(x,t,dm),
\]
and, secondly,
\begin{eqnarray}\nonumber
\int_{\R^d} \int_0^\i J(x,n,t)f(x,t,dn) dx&=&\int_0^\i\int_{\R^d}
h_n(x)J(x,n,0)dxdn  \label{wmsol} \\
&&+\int_0^t\int_0^\i \int_{\R^d}\frac{\partial J}{\partial t}(x,n,s)f(x,s,dn)dxds\nonumber\\
&&+\int_0^t\int_{\R^d} \int_0^\i d(n)\triangle J(x,n,s)f(x,s,dn)dx ds\\ &&+\int_0^t\int_{\R^d} \int_0^\i \int_0^\i \b(m,n)\tilde
J(x,m,n,s)f(x,s,dn)f(x,s,dm)dxds. \nonumber
\end{eqnarray}
\end{definition}
\noindent{\bf Remark:} The requirement  $g_\ell, h_\ell\in L^1_{loc}$
is made in order to guarantee the existence of the integrals in 
(\ref{wmsol}).

We are now ready to state the main result of this article.

\begin{theorem}\label{T:1} Consider the deterministic or random model in some dimension $d \geq 3$. 
 Assume Hypotheses 1.1 and 1.2. 
If $\mathcal{P}$ is any
limit point of $\mathcal{P}^{\epsilon}$, then $\mathcal{P} $ is
concentrated on the space of measures $g (dx,dn,t)dt= f(x,t,dn) dxdt $
which are absolutely continuous with respect to Lebesgue measure $dx\times dt$,
with $f$ solving the system of partial differential equations
\eqref{eq1.1} in the sense of (1.6). The quantities $\b:(0,\infty)
\times (0,\infty) \rightarrow (0,\infty)$ are specified by the
formula
\begin{equation}\nonumber
\b(n,m)=\alpha(n,m) \int_{\mathbb{R}^{d}}V(x)
\left[1+u(x;n,m)\right]dx ,
 \end{equation}
 where, for each pair $(n,m) \in
(0,\infty) \times (0,\infty)$, $u(\cdot)=u(\cdot;n,m):\mathbb{R}^{d} \rightarrow
(0,\infty)$ is the unique solution of
\begin{equation} \label{eq1.12}
\triangle u(x)= \frac{\alpha(n,m)}{d(n)+d(m)}
V(x)\Big[1+u(x) \Big],
\end{equation}
satisfying $u(x) \rightarrow 0$ as $|x| \rightarrow \infty.$
\end{theorem}

\noindent
{\textbf {Remarks.}}
\begin{itemize}
\item
 The continuity with respect to $m$ and $n$ and other important properties of $u(\cdot;n,m)$
will be stated in Lemma 4.2 of Section 4. In particular $u\in[-1,0]$, which implies that $\b>0$
because $u$ is not identically zero. 
It follows from Lemma 4.2 that $\b$ is a continuous function. We also refer to
the last section of \cite{HR1} in which several properties
 of $\b$ are established. In particular, it is shown that 
 $\b\le \a$ and $\beta(n,m)\le {\rm Cap}(K)(d(n)+d(m))$, where $K$ denotes the support of the function $V$ and 
${\rm Cap}(K)$ denotes the Newtonian capacity of the set $K$. (See \cite{HR1} for the definition of Newtonian capacity.)
\item
To simplify our presentation, we assume that all particles have the same ``radius''.
However, in a more realistic model, we may replace $\e^{-2}V( \e^{-1} (x_i-x_j) )$ with 
$\e^{-2}V(\e^{-1} (x_i-x_j) ;m_i,m_j)$, where $V(a;n,m)=(r(n)+r(m))^{-2}V(a/(r(n)+r(m)))$
 and $r(n)$ is interpreted as the radius of a particle 
of mass $n$. 
 Our method of proof applies even when we allow such a radial dependence and we can prove Theorem 1.1 provided that 
$r(n)=n^\chi$ with $\chi<(d-2)^{-1}$ (when $d \geq 3$). In fact, we anticipate that, if $\chi>(d-2)^{-1}$, then, at least in
 the case of a sufficiently large initial condition, the particle densities no longer approximate a solution of (\ref{eq1.1}) 
in which the mass $\int_0^{\infty} \int_{\mathbb{R}^d} m f_m(x,t) dx dt$ is conserved throughout time.
We refer to \cite{R} and the introduction of \cite{HR1}  for a more thorough discussion.
\end{itemize}

Our second result shows that the macroscopic density is absolutely continuous with respect to Lebesgue measure $dn$. 
We will require \\
\noindent{\bf Hypothesis 1.3.}
There exists a continuous function $\tau:(0,\infty) \to (0,\infty)$
for which  $\int_0^\i\tau(n)dn=1$, with 
$$
\int_0^\i\int_{\R^d} \left(|x|^2+|\log \tau(n)|+|\log h_n|\right)h_n\ dxdn<\i.
$$
and
\begin{equation}\label{defrho}
\int_0^\i\int_{\R^d} \rho(n)h_n(x)dxdn<\i,
\end{equation}
where
$$
\rho(n)=\int_0^n\a(m,n-m)\frac{\tau(m)\tau(n-m)}{\tau(n)}dm.
$$
We also assume that $D=\sup_m d(m)<\i$.

\noindent
{\textbf {Remark.}} For a simple example for $\tau$, consider $\tau(n)=(n+1)^{-2}$. If for example $\a(m,n)\le C(m+n)$, then 
$\rho(n)\le Cn$ and \eqref{defrho} requires that the total mass to be finite initially.

\begin{theorem}\label{thmtwo} 
Assume that the model is random, and that the law $\gamma_N$ of the initial total particle number has a Poisson distribution.
 Assume also Hypothesis 1.3. Then every limit point $\cP$ of the sequence
 $\cP^\e$ is concentrated on measures that take the form $g(dx,dn,t)dt
=f_n(x,t)dndxdt$, where $f$ solves (\ref{eq1.1}). Moreover, there exists a
constant $C$, that may be chosen independently of $\cP$, such that 
\begin{equation}
\int_{\cX}\left[\int_0^\i \int_{\R^d} 
 \psi(f_n(x,t))r(x,n)\ dxdn\right]\ \cP(d\mu)\le C,
\end{equation}
for every $t$, where $\psi(f)=f\log f-f+1$ and $r(x,n)=(2\pi)^{-d/2}\exp(-|x|^2/2)\tau(n)$. 
\end{theorem}
\noindent{\bf Remark.}
At the expense of discussing some extra technicalities, the proof of 
Theorem \ref{thmtwo} might include the random model with some other choice of $\gamma_N$. We only need to assume that
 for every positive $\l$, there exists a constant $a(\l)$, such that $\log \mathbb E_N \exp(\l \cN )\le Na(\l)$.  

Theorem \ref{thmtwo} is proved by firstly establishing an entropy bound for the distribution of ${\mathbf q}(t)$,
 and then  using large deviation techniques
to deduce that any limit point $\cP$ of the sequence $\cP^\e$ is concentrated on the space of measures
$g(dx,dn,t)dt=f_n(x,t)dxdndt$. For this, we simply follow the classical work of Guo-Papanicolaou-Varadhan \cite{GPV}.
Even though our result is valid for more general initial randomness, we prefer to state and 
prove our results for Poisson-type distributions, thereby focussing on the main idea of the method of proof.

The function $\tau:(0,\i)\mapsto(0,\i)$ 
appearing in Hypothesis 1.3 is used to define a reference measure with respect to which 
 the corresponding entropy per particle
is uniformly finite as $\e\to 0$.
 For simplicity, we take the reference measure $\nu_N$ which induces a Poisson law of intensity $1$ 
for $\cN$ and whose conditional measure 
$\nu_N(\cdot|\cN(\mathbf q)=k)$ is given by
\begin{equation}
\prod_{i=1}^k r(x_i,m_i)dx_idm_i.
\end{equation}
The entropy per particle is uniformly finite, because
the first part of Hypothesis 1.3 implies that
\[
\sup_N \e^{d-2}\int F^0\log F^0 d\nu_N<\i,
\]
where $F^0(\mathbf q)\nu_N(d\mathbf q)$ denotes the law of $\mathbf q(0)$.
The second part of Hypothesis 1.3 will be used to control the time derivative of the entropy.

We now comment on the possible uniqueness of the solution that the microscopic model approximates. 
We expect to have a unique solution of the system  \eqref{eq1.1} for
the initial condition $h$ as above. However, with the aid of the
arguments of \cite{HR3} and \cite{R2}, we know how to establish this uniqueness
only if we assume  
 that the initial condition satisfies the bound
\begin{equation} \label{eq1.8}
\int_0^{\infty}n^b\|h_n\|_{L^{\infty}}\ dn<\infty,
\end{equation}
for sufficiently large $b=b(a)$ (see \cite{HR3} and \cite{R2} for an expression for $b(a)$). Using this
uniqueness,  we can assert that in fact the limit
$\mathcal{P} $ of $\mathcal{P}^{\epsilon}$ exists and is
concentrated on the single measure $\mu(dx,dn,dt)=f_n(x,t)dxdndt$, where
$f$ is the unique solution to $\eqref{eq1.1}$. As a corollary we
have,
\begin{cor}\label{C:1}
Assume that Hypotheses 1.1, 1.2 and 1.3 hold and that
\eqref{eq1.8} holds for sufficiently large $b$. Let
$J:\mathbb{R}^{d} \times (0,\infty)\times[0,\i)  \rightarrow \mathbb{R}$ be a
bounded continuous function of compact support. Then, 
\begin{equation}\label{eq1.7}
\limsup_{N \rightarrow \infty} \mathbb{E}_{N} \Bigg|\int _{\mathbb{R}^{d}} \int_0^\i \int^{\infty}_{0}  J(x,n,t) \big(\mu(dx,dn,dt)
-f(x,n,t) dx dndt\big)   \Bigg|=0.
\end{equation}
In \eqref{eq1.7}, $f:\mathbb{R}^{d} \times [0,\infty)\times
[0,\infty) \rightarrow [0,\infty) $ denotes the unique solution to
the system \eqref{eq1.1} with  the initial data
$f_{n}(\cdot,0)=h_{n}(\cdot)$.
\end{cor}

The paper contains an appendix that discusses the scalings available in the
Smoluchowski equations in their continuous form. Examining these scalings
produces an heuristic argument for the regime of choices of the asymptotic
behaviour of the input parameters $\beta:(0,\infty)^2 \to (0,\infty)$ and
$d:(0,\infty) \to (0,\infty)$ for which a solution (\ref{eq1.1}) will see
most of the mass depart from any given compact subset of $(0,\infty)$  
as time becomes high.

To outline the remainder of the paper: in Section 2, 
we explain the strategy of the proof, giving an alternative overview to
that presented in \cite{HR1}. In this section, we also show how the microscopic coagulation rate is comparable to the product of 
densities and may be replaced with an
expression that is similar to the term $Q$ in (\ref{eq1.1}) (see Theorem~2.1).
 The main technical step for such a replacement is a regularity property of the
coagulation and is stated as Proposition~2.1. In Section~2 the proof of Proposition~2.1 is reduced to a collection of bounds that are stated
as Lemma 2.1. In Section 3,
we establish the crucial 
 correlation bound (1.4).  In Section 4, the proof of Lemma 2.1 is carried out with the aid of
the correlation bounds of Section 3. In Section 5, we show how the correlation bounds can be used to establish $L^p$-type bounds 
on the macroscopic densities. 
Sections 6 and 7 are devoted to the proofs of Theorems 1.1 and 1.2 respectively.

\noindent{\bf Acknowledgments.} 
We thank James Colliander and Pierre Germain for valuable comments that
relate to the discussion in the appendix. We also thank an anonymous referee for a number of useful suggestions and comments.

\end{section}
\setcounter{equation}{0}

\begin{section}{An outline of the proof of the main theorem}
Our aim in this section is to outline the proof of the principal result,
Theorem 1.1. The overall scheme of the proof is the same as that
presented in \cite{HR1}, and the reader may wish to consult Section $2$ of
that paper for another overview.

Our goal is to show that the empirical measures 
$g^\epsilon(dx,dn,t)$ converge to $f(x,t,dn)dx$, where $f$ is some
measure-valued weak solution of Smoluchowski's equation (\ref{eq1.1}). 
To this end, we choose a smooth test function
$J:\mathbb{R}^{d} \times(0,\infty) \times[0,\infty) \rightarrow
\mathbb{R}$ of compact support and consider the expression
\begin{equation} \nonumber
Y(\mathbf q,t)= \epsilon ^{d-2} \sum_{i \in I_{{\mathbf q}}} J(x_{i},m_{i},t).
\end{equation}
Evidently,
\[
Y({\mathbf q}(t),t)=\int  J(x,n,t)g^\epsilon(dx,dn,t).
\]
Note that
\begin{equation}\label{eq3.1}
Y({\mathbf q}(T),T)=Y({\mathbf q}(0),0)+\int^{T}_{0} \left(\frac{\partial Y}{\partial t}+
\mathbb{A}_{0}(Y) +  \mathbb{A}_{c}(Y)\right)({\mathbf q}(t),t) dt 
+ M_{T},
\end{equation}
where $M_T$ is a martingale, where the free-motion term $\mathbb{A}_{0}Y$ equals
\[
\mathbb{A}_{0}Y({\mathbf q},t) 
= \epsilon^{d-2} \sum_{i \in I_{{\mathbf q}}}
 d(m_{i})\bigtriangleup _{x_{i}} J(x_{i},m_{i},t)=\int  d(n)\bigtriangleup_x J(x,n,t)g^\epsilon(dx,dn,t).
\]
and where the collision term
$\mathbb{A}_{c}Y$ is equal to  
\begin{equation} \label{eq3.2}
\mathbb{A}_{c}Y(\mathbf q,t)    =\epsilon^{d-2} \sum_{i,j \in I_{q}}
\alpha(m_{i},m_{j}) V_{\epsilon}(x_{i}-x_{j}) \hat J(x_i,m_i,x_j,m_j,t),
\end{equation}
with $V_\e(x)=\e^{-2}V(x/\e)$, and $\hat J(x_i,m_i,x_j,m_j,t)$ given by 
\begin{equation}
 \frac{m_{i}}{m_{i}+m_{j}} J(x_{i},m_{i}+m_{j},t)
 +\frac{m_{j}}{m_{i}+m_{j}} J(x_{j},m_{i}+m_{j},t) -J(x_{i},m_{i},t)-J(x_{j},m_{j},t) .
 \end{equation}
Our approach is simply to understand which terms dominate in (\ref{eq3.1})
when the initial particle number $N$ is high, and, in this way, to see that
the equation (\ref{wmsol}) emerges from considering (\ref{eq3.1}) in the high
$N$ limit. Clearly, we expect the last two terms in (\ref{wmsol}),
corresponding to free-motion and collision, to arise from the terms in
(\ref{eq3.1}) in which the operators $\mathbb{A}_0$ or $\mathbb{A}_C$ act. The
time-derivative terms in (\ref{wmsol}) and (\ref{eq3.1}) also naturally
correspond. And indeed, 
the sum of the second and third terms on the right-hand side of (2.1) is already expressed in terms of the empirical
 measure and corresponds to 
the macroscopic expression
\[
\int_0^T\int_0^\infty\int \left(\frac{\partial}{\partial t}+d(n)\bigtriangleup_x\right)J(x,n,t)f(x,t,dn)dxdt.
\] 
As we will see in Section 6, the term martingale $M_{T}$  vanishes as
$\epsilon\to 0$.
The main challenge comes from the fourth term on the right-hand side of
(2.1), the collision term. How does its counterpart in (\ref{wmsol}) emerge
in the limit of high initial particle number? 
To answer this, we need to understand how to express the time-integral of
changes to $Y({\mathbf q},t)$ resulting from all the collisions occurring in the
microscopic model. To do so, it is
natural to introduce the quantity
\[
f^\delta(x,dn;\mathbf q)=\epsilon^{d-2} \sum_{i \in I_{{\mathbf q}}} \delta^{-d}
\xi\Big(\frac{x_{i}-x}{\delta}
\Big)\d_{m_i}(dn),
\]
where $\xi:\mathbb R^d\to [0,\infty)$ is a smooth function of compact support
with $\int_{\R^d} \xi dx=1$.
For $\delta > 0$ fixed and small, $f^\delta$ in essence counts the number
of particles in a small macroscopic region about any given point, this region
having diamater of order $\delta$. To find the analytic collision term in
(\ref{wmsol}) from its microscopic counterpart in (\ref{eq3.1}), we must
approximate the time integral of $\mathbb{A}_cY(\mathbf q(t),t)$ by some functional of the
macroscopically smeared particle count $f^\delta$, in such a way that the
approximation becomes good if we take the smearing parameter $\delta \to 0$
after taking the initial particle number $N$ to be high.
This is achieved by the following important result, in which we write 
$\Gamma(\mathbf q,t)= 
\mathbb{A}_{c}Y(\mathbf q,t)$.
\begin{theorem}\label{T:3} Assume that the function $\hat J(x,m,y,n,t)$
vanishes when $t>T$, or $m+n<L^{-1}$, or $\max({m,n})>L$.
Then
\[
\lim_{\delta\to 0}\limsup_{N\to\infty}\mathbb E_N\left|\int_0^T\left[\Gamma(\mathbf q(t),t)-
\hat\Gamma^\d(\mathbf q(t),t)\right]dt\right|=0,
\]
with
\begin{equation}\label{hatgam}
\hat\Gamma^\d(\mathbf q,t)=\int_{\R^d} \int_{\R^d} \int_0^\i \int_0^\i
\alpha(m,n)U_{m,n}^\epsilon(w_1-w_2)\hat J(w_1,m,w_2,n,t)f^\delta(w_1,dm;\mathbf q)f^\delta(w_2,dn;\mathbf q) dw_1dw_2,
\end{equation}
where we set
\begin{equation}\nonumber
U_{m,n}(x)=V(x)\big[1+u (x;m,n)\big],\ \ \ U_{m,n}^\epsilon(x)=\epsilon^{-d}U_{m,n}(x/\epsilon),
\end{equation}
with $u(\cdot;{m,n})$ being given in Theorem~1.1. 
\end{theorem}

\noindent
{\textbf {Remarks.}}
\begin{itemize}
\item
 Note that even thought $J$ is of compact support, the function $\hat J$ given in (2.3)
is not in general of compact support. 
In fact, if $x_j$ which appears in (2.2) belongs to the bounded support of $J$, then $x_i$ belongs 
to a bounded set because of the presence of the term $V_\e$. The same reasoning does not work for $m_i$ or $m_j$.
Of course if $J(x,n,t)$ vanishes if either $n>L$ or $n<L^{-1}$, then $\hat J(x,m,y,n,t)$ vanishes
if $m+n\le L^{-1}$, or $\max(m,n)>L$. 
However, for Theorem 2.1  we assume that in fact $\hat J$ vanishes even if one of $m$ or $n$ is larger than $L$.
Because of this, we need to show that the contribution of particles with large sizes is small. We leave this issue for Section 6. 
(See Lemma 6.1.)
\item
 As we mentioned in Section~1, the continuity with respect to $m$ and $n$ and other properties of $u(\cdot;m,n)$ will be 
stated in Lemma~4.2.
\end{itemize}

We now explain heuristically why the relation between the cumulative
microscopic coagulation rate $\Gamma(\mathbf q(t),t)$ and its
macroscopically smeared counterpart $\hat\Gamma^\d(\mathbf q(t),t)$ holds.

Here is a naive argument that proposes a form for  $\hat\Gamma^\d(\mathbf
q(t),t)$. In the microscopic model, particles at $(w_1,m)$ and
$(w_2,n)$ are liable to coagulate if their locations differ 
on the scale of $\epsilon$, $\vert w_1 -w_2 \vert =
O(\epsilon)$. If two particles are so located, they coagulate at a Poisson
rate of $\alpha(m,n)V_{\epsilon}(w_1 - w_2)$.
When such a pair does so, it effects a change in $Y(\mathbf q,t)$ of
$\hat{J}(w_1,m,w_2,n)$. The density for the presence of a particle of mass $m$  at location
$w_1$ should be well approximated by the particle count
$f^\delta(w_1,dm)$ computed on a small macroscopic scale. Multiplying
the factors, and integrating over space, we seem to show that the
expression for   $\hat\Gamma^\d(\mathbf q(t),t)$ should be given by
\[
\int_{\R^d} \int_{\R^d} \int_0^\i \int_0^\i
\alpha(m,n)V^\epsilon(w_1-w_2)\hat J(w_1,m,w_2,n,t)f^\delta(w_1,dm;\mathbf q)f^\delta(w_2,dn;\mathbf q)) dw_1dw_2,
\]
where $V^\e(x)=\e^{-d}V(x/\e)$. The  integrand differs from the correct expression in (\ref{hatgam}) by the lack
of a factor
of $1 + \epsilon^{-d} u\big( (w_1 - w_2)/\epsilon  ;m,n \big)$. Why is the preceding argument wrong?
The reason is the following. The joint density for particle presence (of
masses $m$ and $n$) at $w_1$ and $w_2$, (with $\vert w_2 -
w_1 \vert = O(\epsilon)$) is not well-approximated by the product
$f^\delta(w_1,dm)f^\delta(w_2,dn)$, because 
some positive fraction of particle pairs at displacement of order
$\epsilon$
do not in fact contribute, since such pairs were liable to coagulate in the
preceding instants of time, and, had they done so, they would no longer
exist in the model. The correction factor  $1 + \epsilon^{-d} u\big( (w_1 - w_2)/\epsilon  ;m,n \big)$
measures the fraction of pairs of particles, one with diffusion rate
$d(m)$, the other, $d(n)$, that survive without coagulating to reach 
a relative displacement $w_1
- w_2$, and is bounded away from $1$ in a neighbourhood of the origin
of order $\epsilon$.

We note that 
in Theorem 2.1 we have reached our main goals, namely we have produced a quadratic expression of the densities and a function $\alpha U$
which has the macroscopic coagulation propensity $\beta$ for its average. 

The following proposition is the key to proving Theorem \ref{T:3}.
\begin{prop}\label{L:1} Choose $T$ large enough so that $\hat J(\cdot,t)=0$ when $t\ge T$. We have
\begin{equation} \label{eq3.6}
\lim_{|z|\rightarrow 0}\limsup _{\epsilon \downarrow 0}
\mathbb{E}_{N} \left|\int_{0}^{T} \left[\Gamma({\mathbf q}(t),t)-
    \bar{\Gamma}_z({\mathbf q}(t),t)\right]dt\right|=0, 
\end{equation}
where
\begin{equation} \label{eq3.5}
\bar{\Gamma}_z(\mathbf q,t)=\epsilon^{2(d-2)} \sum_{i,j \in I_{\mathbf q}} \alpha(m_{i},m_{j})
U_{m_{i},m_{j}}^\e\left({x_{i}-x_{j}+z}\right)\hat J(x_{i},m_{i},x_{j},m_{j},t).
\end{equation} 
\end{prop}
In the statement, $z$ plays the role of a small macroscopic displacement,
taken to zero after the limit of high initial particle number is taken in
the microscopic model. The proposition shows that the cumulative influence
of coagulations in space and time on $Y(\mathbf q(t),t)$ is similar to that
computed by instead considering pairs of particles at the fixed small macroscopic
distance $z$, with a modification in the coagulation propensity in the
expression (\ref{eq3.5}) being made for the reason just described. 

It is not hard to deduce Theorem 2.1 from Proposition \ref{L:1}.
 We refer to Section 3.5 of \cite{HR1} for a proof of Theorem 2.1 assuming
Proposition \ref{L:1}. See also [10] for a repetition of this proof and more heuristic discussions about the strategy of the proof.

We will prove Proposition \ref{L:1} in the following way.
Define  
\begin{equation*}
X_{z}(\mathbf{q},t)=\epsilon^{2(d-2)} \sum_{i,j\in I_{\mathbf q}}
u^{\epsilon}(x_{i}-x_{j}+z;m_i,m_j)\,\hat J(x_
{i},m_{i},x_{j},m_{j},t),
 \end{equation*}
where $u^{\epsilon}(x;m,n)=\e^{2-d}u(x/\e;m,n)$. Note that $u^\e(x)=u^\e(x;m,n)$ solves
\begin{equation}\label{equpde}
(d(m)+d(n))\Delta u^\e=\a(m,n)(V_\e u^\e+V^\e),
\end{equation}
with 
\[
V_\e(x)=\e^{-2}V(x/\e),\ \ \ \ \ V^\e(x)=\e^{-d}V(x/\e).
\]
The process $\big\{ \big( X_{z} - X_{0} \big)(\mathbf q (t),t): t \geq 0 \big\}$
satisfies
\begin{eqnarray}\label{eq3.6a}
    \big( X_{z} - X_{0} \big) \big(\mathbf q (T),T  \big) =& &
    \big( X_{z} - X_{0} \big) \big( \mathbf q (0),0 \big) +
\int_{0}^{T}{ \Big( \frac{\partial}{\partial t} + \mathbb{A}_{0}
\Big)
\big(X_{z} - X_{0}\big) (\mathbf q (t),t)dt} \\
    &&+ \,
\int_{0}^{T}{ \mathbb{A}_{c} (X_{z} - X_{0}) (\mathbf q (t),t)dt} \, + \, M(T),
\nonumber
\end{eqnarray}
with $\big\{M(t)\,:\,t \geq 0 \big\}$ being a martingale. 
We will see that the form (\ref{eq3.6}) emerges from the dominant terms in
(\ref{eq3.6a}), those that remain after the limit of high initial particle number
$N \to \infty$ is taken. To see this,  
we label the various terms which appear on the right-hand side of (\ref{eq3.6a}).
Firstly, those terms arising from the action of the
diffusion operator: 
\begin{equation*}
 \Big( \frac{\partial}{\partial t} + \mathbb{A}_0 \Big) (X_z -
X_0) = H_{11} + H_{12} + H_{13} +H_{14} +
H_2 + H_3 + H_4,
\end{equation*}
 with
\begin{eqnarray} \nonumber
H_{11}  (\mathbf q,t) &= & \epsilon^{2(d-2)} \sum_{i,j \in
I_{{\mathbf q}}}{\alpha(m_{i},m_{j})
    \Big[ V^{\epsilon} \big( x_{i} - x_{j} + z  \big) -
V^{\epsilon} \big(   x_{i} - x_{j}  \big) \Big]} \hat J(x_{i},m_{i},x_{j},m_{j},t) ,\\
H_{12}(\mathbf q,t)  & = & - \epsilon^{2(d-2)} \sum_{i,j \in
I_{{\mathbf q}}}{\alpha(m_{i},m_{j})
    V_{\epsilon} \big( x_{i} - x_{j}  \big) u^{\epsilon} \big( x_{i} -
    x_{j} ;m_i,m_j
    \big)  \hat J(x_{i},m_{i},x_{j},m_{j},t) } \nonumber ,\\
H_{13}(\mathbf q,t)  & = &  \epsilon^{2(d-2)} \sum_{i,j \in
I_{\mathbf q}}{\alpha(m_{i},m_{j})
     V_{\epsilon} \big( x_{i} - x_{j}+ z \big) u^{\epsilon} \big( x_{i} - x_
{j}+z;m_i,m_j
    \big)  \hat J(x_{i},m_{i},x_{j},m_{j},t) }, \nonumber\\
 H_{14}(\mathbf q,t)  & = &  \epsilon^{2(d-2)} \sum_{i,j \in I_{{\mathbf q}}} d(m_{i})
          \Big[  u^{\epsilon} (x_{i} - x_{j} + z;m_i,m_j ) - u^{\epsilon}
(x_{i} - x_{j};m_i,m_j ) \Big]
     \hat J_{t}(x_{i},m_{i},x_j,m_j,t), \nonumber 
\end{eqnarray}
along with
\begin{eqnarray} \nonumber
H_{2}(\mathbf q,t)  & =  & 2 \epsilon^{2(d-2)} \sum_{i,j \in I_{{\mathbf q}}} d(m_{i})
          \Big[  u^{\epsilon}_{x} (x_{i} - x_{j} + z;m_i,m_j ) - u^{\epsilon}_{x}
(x_{i} - x_{j} ;m_i,m_j ) \Big]
    \cdot \hat J_{x}(x_{i},m_{i},x_j,m_j,t), \\
H_{3}(\mathbf q,t)   & = & -  2 \epsilon^{2(d-2)} \sum_{i,j \in I_{{\mathbf q}}}d(m_{j})
      \Big[  u^{\epsilon}_{x} (x_{i} - x_{j} + z;m_i,m_j  ) -
     u^{\epsilon}_{x} (x_{i} - x_{j};m_i,m_j  ) \Big]
    \cdot \hat{J}_{y} (x_i,m_i,x_{j},m_{j},t) , \nonumber
\end{eqnarray}
and
\begin{eqnarray} \nonumber
H_{4}(\mathbf q,t) &=&\epsilon^{2(d-2)}\sum_{i,j \in I_{{\mathbf q}}}\Big[
u^{\epsilon}(x_{i}-x_{j}+z;m_i,m_j ) - u^{\epsilon}(x_{i}-x_{j};m_i,m_j ) \Big]\\
\nonumber
 && \qquad \qquad \quad \Big[d(m_{i}) \Delta_{x}\hat J(x_{i},m_{i},x_j,m_j,t)
 + d(m_{j})
 \Delta_{y}\hat{J}(x_i,m_i,x_{j},m_{j},t)\Big], \\ \nonumber
\end{eqnarray}
where $\hat J_{x}$ denotes the gradient of $\hat J$ with respect to its first spatial argument,
$\hat J_{y}$  the gradient of $\hat J$ with respect to its second spatial argument, and $\cdot$ the scalar
product. As for those terms arising from the action of the collision
operator,
\begin{equation} \nonumber
\mathbb{A}_{c}(X_{z}-X_{0})(\mathbf{q},t)=
G_{z}^1( \mathbf q ,t)+G_{z}^2(\mathbf q ,t)-G_{0}^1(\mathbf q ,t)-G_{0}^2(\mathbf q ,t),
\end{equation}
where $G_{z}^1 (\mathbf q ,t)$ is set equal to
\begin{eqnarray}\nonumber
&& \sum_{k,\ell\, \in \,I_{{\mathbf q}}}{\alpha(m_{k},m_{\ell}) V_{\epsilon}(
x_{k} - x_{\ell})} \epsilon^{2(d-2)} \sum_{i \,\in \, I_{{\mathbf q}}} \\
\nonumber    && \quad \bigg\{
    \frac{m_{k}}{m_{k} + m_{\ell}} \Big[  u^{\epsilon} ( x_{k} -x_{i}
    + z;m_k+m_\ell,m_i) \hat J( x_{k} , m_{k} + m_{\ell},    x_{i},m_{i},t)  \\ \nonumber
    && \qquad \qquad \qquad  + \, u^{\epsilon} (x_{i} - x_{k}
    + z;m_i,m_k+m_\ell) \hat J( x_{i} , m_{i} ,x_{k},m_{k} + m_{\ell},t)  \Big]
     \\  \nonumber
    && \quad \, \, + \, \frac{m_{\ell}}{m_{k} + m_{\ell}}  \Big[ u^
{\epsilon} ( x_{\ell} - x_{i}
    + z;m_k+m_\ell,m_i) \hat J( x_{\ell} , m_{k} + m_{\ell},x_{i},m_{i},t)  \\ \nonumber
    &&  \qquad \qquad \qquad  + \, u^{\epsilon} ( x_{i} - x_{\ell}
    + z;m_i,m_k+m_\ell) {\hat J}( x_{i} , m_{i} , x_{\ell},m_{k} + m_{\ell},t) \Big]
      \\ \nonumber
    && \qquad  - \ \
    \Big[  u^{\epsilon} (x_{k} - x_{i} + z;m_k,m_i) {\hat J}(x_{k} , m_{k},
    x_{i},m_{i},t) \nonumber \\ \nonumber
    &&  \qquad \qquad \qquad + \,  u^{\epsilon} (x_{i} - x_{k} + z;m_i,m_k)
    {\hat J}( x_{i} , m_{i} , x_{k},m_{k},t) \Big]
      \\ \nonumber
    && \qquad - \ \
     \Big[  u^{\epsilon}(x_{\ell} - x_{i} + z;m_\ell,m_i) {\hat J}( x_{\ell} , m_{\ell},
    x_{i},m_{i},t) \\ \nonumber
    && \qquad \qquad \qquad  + \,  u^{\epsilon} (x_{i} - x_{\ell} + z;m_i,m_\ell) 
{\hat J}( x_{i} ,    m_{i} , x_{\ell},m_{\ell},t) \Big] \bigg\}
      ,
\end{eqnarray}
and where
\begin{equation} \nonumber
G_{z}^2 (\mathbf q ,t)  =  - \epsilon^{2(d-2)} \sum_{k,\ell \,\in I_{{\mathbf q}}}{
\alpha(m_{k},m_{\ell}) V_{\epsilon}(x_{k} - x_{\ell})} u^{\epsilon}
(x_{k} - x_{\ell} + z ;m_k,m_\ell) {\hat J}(x_{k},m_{k},
x_{\ell},m_{\ell},t).
\end{equation}
 The terms in $G_{z}^1$ arise from the changes in the functional
$X_{z}$ when a collision occurs due to the influence of the
appearance and disppearance of particles on other particles that are
not directly involved. Those in $G_{z}^2$ are due to the absence
after collision of the summand in $X_{z}$ indexed by the colliding
particles.

As we take a high $N$ limit in (\ref{eq3.6a}), note that the quantity  
$$
\int_0^T
\Gamma(\mathbf q(t),t) dt =  \epsilon^{2(d-2)} \sum_{i,j \in
I_{{\mathbf q}}}{\alpha(m_{i},m_{j})   
V^{\epsilon} \big(   x_{i} - x_{j}  \big) } \hat J(x_{i},m_{i},x_{j},m_{j},t)
$$
appears, with a negative sign, in the term $H_{11}$. The term $H_{12}$ also
remains of unit order in the high $N$ limit, and would disrupt our aim of
approximating 
$\int_0^T
\Gamma(\mathbf q(t),t) dt$ by $z$-displayed expressions. However, our
definition of $u^{\epsilon}$ (see (\ref{equpde})) ensures that
\begin{equation*}
    H_{12} - G_{0}^2  = 0,
\end{equation*}
so that this unwanted term disappears. The definition of $u^\epsilon$ was
made in order to achieve this. The other term of unit order remaining in
the high $N$ limit is the $z$-displaced $H_{13}$.
Rearranging (\ref{eq3.6a}), we obtain
\begin{eqnarray}\label{estrhs}
\Big| \int_{0}^{T}\,H_{11}(\mathbf{q}(t),t)dt+ \int_{0}^{T}\, H_{13}(\mathbf{q}(t),t)dt \Big|  
&   \leq &|X_{z}-X_{0}|\big({\mathbf q}(T),T\big) +
|X_{z}-X_{0}|\big({\mathbf q}(0),0\big) \nonumber\\  
&& +\int_{0}^{T}\,
\left( |H_{14}|+|H_{2}|+|H_{3}|+ |H_{4}|\right)(\mathbf q(t),t)dt \\ 
 &&+\int_{0}^{T}\, |G_{z}^1-G_{0}^1| (\mathbf q(t),t)dt+\,\int_{0}^{T}\,
|G_{z}^2|(\mathbf q(t),t)dt\nonumber\\
&&+\,\big|M(T)\big|.\nonumber
 \end{eqnarray}
We have succeeded in writing 
$\overline{\Gamma}_z-\Gamma$ in the form $H_{11} + H_{13}$, so that, for Proposition
\ref{L:1}, it remains to 
 prove that the right-hand-side of (\ref{estrhs}) is small enough. 
Firstly,
recall that, by our assumption, the function $\hat J$ is of compact support. We now choose $T$ sufficiently large 
so that $\hat J(x,m,y,n,T)=0$.
As a result, the first term on the right-hand side  vanishes. 
The other bounds we require are now stated.
\begin{lemma}\label{lemL:2}
There exists a constant $C_2=C_2(\hat J,T)$ such that,
\begin{eqnarray} \label{eq3.10}
 && \int_0^T \mathbb{E}_{N}\left(|H_{2}|+|H_{3}|\right) (\mathbf q(t),t) dt\leq C_2|z|^{\frac{1}{d+1}}, \\ 
\label{eq3.11}
&& \int_0^T \mathbb{E}_{N}\left(|H_{4}|+|H_{14}| \right)(\mathbf q(t),t) dt\leq C_2|z|^{\frac{2}{d+1}}, \\ \label{eq3.12}
&& \int_{0}^{T}\, \mathbb{E}_{N}|G_{z}^1-G_{0}^1|(\mathbf q(t),t) dt\leq C_2|z|^{\frac{2}{d+1}},\\ \label{eq3.13}
&& \int_0^T  \mathbb{E}|G_{z}^2|(\mathbf q(t),t) dt\leq C_2\Big(
\frac{\epsilon}{|z|}\Big)^{d-2}, \\ 
\label{eq3.14}
&&\mathbb{E}_{N}|X_{z}-X_{0}|(\mathbf q(0))  \leq C_2 |z|,\\ \label{eq3.15}
&&\mathbb{E}_{N}\big[M(T)^{2}\big]  \leq C_2\epsilon^{d-2}.
\end{eqnarray}
\end{lemma}
These bounds are furnished by the correlation inequality Theorem \ref{T:2}
that is the main innovation of this paper, to whose proof we now turn. 
\end{section}

\setcounter{equation}{0}
\begin{section}{Correlation Bounds}

This section is devoted to the proof of the correlation bound which appeared as (1.4) when $d(\cdot)$ is non-increasing and
takes the form (3.1) in general. Recall the function $\phi$ which appeared in Hypothesis 1.1.
The main result of this section is Theorem~\ref{T:2}.
\begin{theorem}\label{T:2}
For every non-negative bounded continuous function $K:(
\mathbb{R}^{d} )^{k} \,\rightarrow \mathbb{R}$,
\begin{eqnarray} \label{eq2.1}
 &&\mathbb{E}_N \int_0^\i\, \sum_{i_{1},\ldots,i_{k}\in
I_{\mathbf{q}(t)}} \,
K \big(x_{i_{1}}(t),\ldots,x_{i_{k}}(t) \big)\prod_{r=1}^{k}\g_k\big(m_{i_{r}}(t)  
\big) dt\\
 & &\ \ \ \leq  \mathbb{E}_N\,
\sum_{i_{1},\ldots,i_{k}\in
I_{\mathbf{q}(0)}}( {\L}^{m_{i_1}(0),\dots,m_{i_k}(0)} K)\big(x_{i_{1}}(0),\ldots,x_{i_{k}}(0)
\big) \prod_{r=1}^{k}\g_k\big(m_{i_r}(0)\big) , \nonumber
\end{eqnarray}
where all summations are over
distinct indices $i_{1}, \ldots , i_{k} $, the function $\g_k(m)=md(m)^{d/2}\phi(m)^{\frac {kd}2-1}$, 
and the operator ${\L}$ is defined by
\begin{equation}
{\L}^{n_1,\dots,n_k} K(y_1,\dots,y_k)=c_0(kd)\int \prod_{r=1}^kd(n_r)^{-d/2} 
\left(\frac{|y_1-z_1|^2}{d(n_1)}+\dots+\frac{|y_k-z_k|^2} {d(n_k)}\right)^{1-\frac{kd}2}
K(z_1,\dots,z_k) dz_r,
\end{equation} 
where $c_0(kd)=(kd-2)^{-1}\o_{kd}^{-1}$, with $\o_{kd}$ denoting the surface area of the unit sphere in $\bR^{kd}$.
\end{theorem}

Let us make a comment about the form of (\ref{eq2.1}) before embarking on its proof.
Observe that if there were no coagulation, then it would have been straightforward to 
bound the left-hand side of (\ref{eq2.1}) with the aid of the diffusion semigroup even if we allow a
 function $K$ that depends on the
 masses of 
particles. Indeed, if $S_t^{m_{i_1},\dots,m_{i_k}}$ denotes the diffusion semigroup associated with 
particles $(x_{i_1},m_{i_1}),\dots,(x_{i_k},m_{i_k})$, then $\int_0^\i S_t^{m_{i_1},\dots,m_{i_k}}dt$ 
is exactly the operator $\L^{m_{i_1},\dots,m_{i_k}}$.
What (3.1) asserts is that a similar bound is valid in spite of coagulation provided that we allow only a very special
dependence on the masses of particles.

\noindent{\bf Proof of Theorem 3.1.} Let us define 
\[
G(\mathbf q)=\sum_{i_{1},\ldots,i_{k}\in
I_{\mathbf{q}}}({\L}^{m_{i_1},\dots,m_{i_k}}K)\big(x_{i_{1}},\ldots,x_{i_{k}}
\big) \prod_{r=1}^{k}\g_k(m_{i_{r}}) .
\]
Recall that the process
$\mathbf{q}(t) $
is a Markov process with generator
$\mathbb{L}=\mathbb{A}_{0}+\mathbb{A}_{c}$ where
$\mathbb{A}_{0}=\sum_{i\in {I_{\mathbf q}}}d(m_{i})\Delta_{x_{i}}.$ By Semigroup Theory,
\begin{equation}
\label{eq2.3}
\mathbb{E}_NG\big(\mathbf{q}(t)\big)=\mathbb{E}_NG\big(\mathbf{q}(0)\big)
+\mathbb{E}_N\int_{0}^{T} \mathbb{L}G(\mathbf{q}(t))dt.
\end{equation}
We have
\begin{equation}
\cA_0G(\mathbf q)=-\sum_{i_{1},\ldots,i_{k}\in
I_{\mathbf{q}}}K\big(x_{i_{1}},\ldots,x_{i_{k}}
\big) \prod_{r=1}^{k}\g_k(m_{i_{r}}) .
\end{equation}
This and the assumption $K\ge 0$ would imply (\ref{eq2.1}) provided that we can show
\begin{equation}
\cA_cG\le 0.
\end{equation}

To prove (3.5), let us study the
effect of a coagulation between the $i$-th and $j$-th particle on $G$.
We need to study three cases separately:
\begin{itemize}
\item
 $i,j
\notin \{i_{1},\ldots,i_{k}\},$
\item
  $i,j \in
\{i_{1},\ldots,i_{k}\},$ 
\item
  only one of  $i$ or $j$ belongs to $ \{i_{1},\ldots,i_{k}\}$.
\end{itemize}

 If the first case occurs, then $(i,j)$-coagulation does not
affect the term indexed by  $\{i_{1},\ldots,i_{k}\}$   in $G(\mathbf q)$.

If the second case occurs,
then we need to remove those terms in the summation for which
$\{i,j\}\subseteq \{i_{1},\ldots,i_{k}\}$. This contributes
negatively to $\mathbb{A}_{c}G(\mathbf q)$,
because $K\geq 0$. This total contribution for this case is given by
\begin{eqnarray}\nonumber
 & & -\sum_{i,j
\in{I_{\mathbf{q}}}}
V_{\epsilon}(x_{i}-x_{j})\, \alpha(m_{i},m_{j})\\
\nonumber &&\ \ \cdot \sum_{i_{1},\ldots,i_{k}} 1 \!\! 1\big( i,j
\in \{i_{1},\ldots,i_{k} \} \big)
\big( \L^{m_{i_1},\dots,m_{i_k}} K
\big)(x_{i_{1}},\ldots,x_{i_{k}}) \prod_{r=1}^{k} \g_k(m_{i_r}) . 
\end{eqnarray}

If the third case occurs, then only
one of $i,j$ belongs to $\{i_{1},\ldots,i_{k}\}$. For example, either  $i=i_{1}$, and
$j \notin \{i_{1},\ldots,i_{k}\}$, or
 $j=i_{1}$, and
$i \notin \{i_{1},\ldots,i_{k}\}$. In this case, the contribution is still
non-positive because after the coagulation the expression
\begin{eqnarray}
\nonumber Y_1 &=& \sum_{i_{2},\ldots,i_{k}}
\big(\L^{m_i,m_{i_2},\dots,m_{i_k}} K
\big)(x_{i},x_{i_{2}},\ldots,x_{i_{k}})  \g_k(m_{i})   \prod_{r=2}^{k} \g_k(m_{i_r})   \\
&& \quad +\sum_{i_{2},\ldots,i_{k}}\big(\L^{m_j,m_{i_2},\dots,m_{i_k}} K
\big)(x_{j},x_{i_{2}},\ldots,x_{i_{k}}) \g_k(m_{j})   \prod_{r=2}^{k} \g_k(m_{i_r}) ,  
\nonumber
\end{eqnarray}
is replaced with the expression $Y_2$ which is given by
\begin{eqnarray}\nonumber
&& \frac{m_{i}}{m_{i}+m_{j}}\sum_{i_{2},\ldots,i_{k}}\big(\L^{m_{i}+m_{j},m_{i_2},\dots,m_{i_k}} )K
\big)(x_{i},x_{i_{2}},\ldots,x_{i_{k}}) \g_k(m_{i}+m_j)   \prod_{r=2}^{k} \g_k(m_{i_r})  \\ 
&&+\frac{m_{j}}{m_{i}+m_{j}}\sum_{i_{2},\ldots,i_{k}}\big(\L^{m_i+m_j,m_{i_2},\dots,m_{i_k}} )K
\big)(x_{j},x_{i_{2}},\ldots,x_{i_{k}}) \g_k(m_i+m_{j})   \prod_{r=2}^{k} \g_k(m_{i_r}) .  
 \nonumber
\end{eqnarray}
For (3.5), it suffices to show that $Y_2\le Y_1$. For this,  it suffices to show that for every positive $m,n,A$ and $B$,
\begin{equation}
\phi(m+n)^{\frac{kd}2-1} \left[A\frac{d(m)}{d(m+n)}+B\right]^{1-\frac{kd}2}\le \phi(m)^{\frac{kd}2-1} [A+B]^{1-\frac{kd}2}.
\end{equation}
We are done because the assertion (3.6) for fixed $m,n$ and all positive $A$ and $B$
is equivalent to 
the inequalities 
$$
\phi(m) d(m) \geq \phi\big( m + n \big) d \big( m + n \big),
$$
and
$$
\phi(m) \geq \phi (m+n),
$$
both being satisfied, and these are 
true for all choices of $m$ and $n$ by Hypothesis 1.1.
\qed

\begin{cor}\label{C:2}
For every non-negative bounded continuous
function  $K$,
\begin{eqnarray} \label{eq1.13}
 &&\e^{k(d-2)}\mathbb{E}_N \, \int_0^T\sum_{i_{1},\ldots,i_{k}\in
I_{\mathbf{q}(t)}} \,
K \big(x_{i_{1}}(t),\ldots,x_{i_{k}}(t) \big)\prod_{r=1}^{k}\g_k\big(m_{i_{r}}(t)\big) dt\\
 & &\ \ \ \leq  c_0(kd) \int K(x_1,\dots,x_k) \prod_{r=1}^k 
\left(\bar h_k*\l_k\right )(x_r) dx_r,\nonumber
\end{eqnarray}
where $\bar h_k=\int_0^\i n\phi(n)^{\frac{kd}2-1}d(n)^{\frac d2-\frac 1k}h_n\ dn$ and
$\l_k (w) =|w|^{\frac 2k-d}.$
\end{cor}
\textbf{Proof. } 
From  the elementary inequality $a_1\dots a_k\le (a_1^2+\dots+a_k^2)^{k/2}$, we deduce that the kernel $\l^{n_1,\dots,n_k}$ of the operator
$\L^{n_1,\dots,n_k}$ is bounded above by
\[
\l^{n_1,\dots,n_k}(z_1,\dots,z_k)\le c_0(kd) \prod_{r=1}^k|z_r|^{\frac 2k-d}d(n_r)^{-\frac 1k} .
\]
This and (3.1) imply (3.7).
\qed

We end this section with two lemmas concerning the first condition in Hypothesis 1.1.
\begin{lemma}
Suppose the function $d(\cdot)$ has a finite negative variation in an interval $[a,b]\subset (0,\i)$. Then 
there exists a positive continuous function $\phi$ such that $\phi$ and $\phi d$ are non-increasing in the interval $[a,b]$.
\end{lemma}
\noindent{\bf{Proof. Step 1.}} Firstly, we assume that there exist points $a_0=b>a_1>\dots>a_{\ell-1}>a_{\ell}=a$ such that $d(\cdot)$ is 
monotone on each interval $[a_{i},a_{i-1}]$, $i=1,\dots,\ell$. For the sake of definiteness, let us assume that $d(\cdot)$ is
non-decreasing (non-increasing) in $[a_i,a_{i-1}]$, if $i$ is odd (even). In this case, we can construct a continuous $\phi$
as follows: 
Define $A_0 = A$ and $A_k = A\prod_{i=1}^{k} \frac{d(a_{2i})}{d(a_{2i - 1} )}$
for $k \ge 1$.
For $x \in\big[ a_{2k+1} , a_{2k} \big]$ and $k \ge 0$,
we set $\phi(x) = \frac{A_k}{d(x)}$.
For $x \in \big[ a_{2k} , a_{2k - 1} \big]$ and $k \ge 1$,
we set $\phi(x) = \frac{A_{k-1}}{d(a_{2k-1})}$. 

\noindent{\bf Step 2.} Let $d$ be a continuous positive function. Approximate $d$ in $L^{\infty}$ by a sequence of continuous 
piecewise monotone functions
$\{d_n\}$. To simplify the presentation, we assume  that each $d_n$ is as in Step 1. That is, $d_n$ increases near the end point $b$. 
 Let us write $\phi_n$ for the corresponding $\phi$, and let $c_n$ denote the number of 
intervals in the partition (so that $a_{c_n} = a$). It remains to show that
the sequence $\{\phi_n\}$ has a convergent subsequence. 
Since each $\phi_n$ is non-increasing, we may appeal to the Helley Selection Theorem. For this we need to make sure that the sequence 
 $\{\phi_n\}$ is bounded. 
Note that $\sup_{x \in [a,b]} \phi_n(x) = \phi_n(a) = \phi_n\big(a_{c_n}\big)$.
Set $D_n = A_{\frac{c_n - 1}{2}}$ if $c_n$ is odd and
$D_n = A_{\frac{c_n}{2} - 1}$ if $c_n$ is even. We readily see that $\phi_n \big( c_n \big)
 \leq \big( \inf_{x \in [a,b]} d(x) \big)^{-1} D_n$,
 whatever the parity of $c_n$. The infimum being positive, we require that $\sup_{n \in \mathbb{N}}{D_n} < \infty$. 
For any $k \in \mathbb{N}$ for which $A_k$ is defined, 
we may take the logarithm of $A_k$ to produce a sum and observe that $d(\cdot)$ is non-increasing on the intervals
$[a_{2i},a_{2i-1}]$. Hence, $\log A_k$ measures the negative variation of  the function $\log d$ on the interval $\big[ a_{2k}, b \big]$.  
 Since $d$ is uniformly positive, $\sup_n D_n < \infty$ is implied by 
the function $d$ having a finite negative variation.
\qed 

\begin{lemma}
Suppose the function $\log d(\cdot)$ has a finite negative variation in an interval $[n_0,\i)$ with $n_0>0$. Then 
there exists a function positive continuous $\phi$ such that $\phi$ and $\phi d$ are non-increasing in the interval $[n_0,\i)$.
\end{lemma}
\noindent{\bf Proof.} The proof is very similar to the proof of Lemma 3.1. First we assume that $d$ is piecewise
monotone. This time we set $\phi(n_0)=A$ and define $\phi$ continuously so that $\phi$ is constant when $d$ decreases and 
$\phi$ is a constant multiple of $d^{-1}$ when $d$ increases. Since $\phi$ is non-increasing, we may end with a function which crosses 
0 and becomes negative. This can be fixed by adjusting $A=\phi(n_0)$, only if $\phi$ is bounded below. As in the proof of Lemma 3.1,
we can readily see that $\phi$ is bounded below if the total negative variation of $\log d$ is finite.
\qed

 Note that in the statement of Lemma 3.2 we can not drop $\log$ because on the infinite interval $[n_i,\i)$ the function $d(\cdot)$ could
take arbitrarily small values.

\end{section}
\setcounter{equation}{0}
\begin{section}{Proof of Lemma 2.1}
The strategy of the proof of Lemma 2.1 is the same as the one used to prove the analogous inequalities in
\cite{HR1}. The only difference is that we only need to use our correlation bound Corollary 3.1 to get the bounds (2.10--15). 
For (2.10) and (2.11) we need to apply Corollary 3.1 for $k=2$.  Corollary 3.1 in the case $k=3$ will be used for (2.12).
 As for (2.15)
all cases $k=2,3,4$ will be employed. 
We omit the proof of the inequalities  (2.13) and (2.14) because they can be established by a verbatim 
argument as in \cite{HR1}. In fact the proof (2.14) is straightforward because we are dealing with a calculation involving the initial
configuration. For this, however, a suitable bound on the function $u^\e$ would be
 needed that will be stated as a part of Lemma 4.2 below.
The same bound and Lemma 4.1 below will imply (2.13). 

The main ingredients for the proof of inequalities (2.10) and (2.11) are Corollary 3.1 (with $k=2$), certain bounds 
on $u^\e$ and $u_x^\e$ (which will appear in Lemma 4.2), and Lemma 4.1 below.
The straightforward proof of Lemma 4.1 is also omited and  can be proved in exactly the same way we proved 
Lemma 3.1 of \cite{HR1}. 
\begin{lemma}\label{lemma4.1} For any $T \in [0, \infty)$,
\begin{equation}\nonumber
\mathbb{E}_{N}\, \int^{T}_{0}\, \epsilon^{d-2} \sum_{i,j \in
I_{{\mathbf q}}(t)} \, \alpha (m_{i}(t),m_{j}(t))V_{\epsilon}(x_{i}(t)-x_{j}(t))dt \leq Z.
\end{equation}
\end{lemma}

As for the remaining inequalities, we only establish (2.12) and (2.15) because these are the most technically
 involved cases and the same
 idea of proof applies to (2.10) and (2.11).

We now state our  lemma about the functions $u$ and $u^\e$. Recall that
 $u^{\epsilon}(x;n,m)=\epsilon^{2-d}u(x/\epsilon;n,m)$
 where $u$ satisfies
 \begin{equation} \nonumber
 \triangle
 u(x;n,m)=\a'(n,m) V(x)\big[1+u(x;n,m)
 \big],
 \end{equation}
 with $u(x;n,m) \rightarrow 0$ as $|x| \rightarrow \infty$, and
\[
\a'(n,m):=\frac{\alpha(n,m)}{d(n)+d(m)}.
\]
For our purposes, let us write $w^a$ for the unique solution of 
\begin{equation} \nonumber
 \triangle
 w^a(x)=a V(x)\big[1+w^a(x)
 \big],
 \end{equation}
with $w^a(x) \rightarrow 0$ as $|x| \rightarrow \infty$. Of course, if we choose $a=\a'(n,m)$, then we obtain $u(x;n,m)$.
 We choose the constant $C_{0}$ so that $V(x)=0$ whenever $|x| \geq
 C_{0}$. 

 \begin{lemma} There exists a constant $C_3$
 for which the following bounds hold. 
\begin{itemize}
\item $-1\le w^a(x)\le 0$ and for $x \in \mathbb{R}^{d}$,
\begin{eqnarray}\nonumber
|w^a(x)| &\leq&
        C_3a\min\{|x|^{2-d},1\},\\
\nonumber
|w^a_x(x)| &\leq&
        {C_3}a\min\{|x|^{1-d},1\}.
\end{eqnarray}
\item for $x \in \mathbb{R}^{d}$ satisfying $|x| \geq \max
\big\{2|z|+C_{0}\epsilon , 2C_{0}\epsilon \big\}$,
\begin{equation} \label{eq4.3}
\big |u^{\epsilon}(x+z;n,m)-u^{\epsilon}(x;n,m)\big | \leq
{C_3\a'(n,m)|z|}{|x|^{1-d}}
\end{equation}
and
\begin{equation} \label{eq4.4}
\big|u^{\epsilon}_x(x+z;n,m)- u^{\epsilon}_x(x;n,m)\big
| \leq {C_3\a'(n,m)|z|}{|x|^{-d}}.
\end{equation}
\item the function $w^a$ is  differentiable with respect to $a$ and  
$a^{-1}w^a\le \frac{\partial w^a}{\partial a}\le 0$.
\end{itemize}
\end{lemma}

\noindent
\textbf{Proof.} The proof of the first and second parts
 can be found in Section 3.2 of \cite{HR1} and we do not repeat it here. As for the third part, recall that 
the function $w^a$ is uniquely determined by the equation
\begin{equation}\label{eq4.31}
w^a(x)= -c_0 a\int_{\mathbb{R}^d} |x-y|^{2-d}V(y)(1+w^a(y))dy,
\end{equation}
where $c_0=c_0(d)=(d-2)^{-1}\o_d^{-1}$, with $\o_d$ denoting the surface area of
the unit sphere $S^{d-1}$. We wish to show the regularity of the function $w^a$ with respect to the variable $a$.
In fact the existence of the unique solution to (4.3) was established in \cite{HR1} using the
Fredholm Alternative Theorem. 
To explain this, let us pick  a
bounded continuous
function $R$ such that  $R > 0$, with
\[
\int_{\mathbb{R}^d} R(x)dx = \i,\ \ \ \int_{|x| \ge 1}
R(x)|x|^{4-2d}dx <
\i.  
\]
Define
\begin{equation}
\nonumber
\cH = \left\{ u : \mathbb{R}^d \to \mathbb{R}: \mbox{ $u$ is measurable and }\int_{\mathbb{R}^d} 
u^2(x) R(x) dx <
\i\right\}.
\end{equation}
Observe that $\cH$ is a Hilbert space with respect to
the inner product
\[
\<u,v\> = \int_{\mathbb{R}^d} u(x)v(x) R(x)dx.
\]
Note that if $w^a$ solves (4.3), then, defining $\cF:\cH\mapsto \cH$ by
$$
\cF (\omega) = c_0 \int \vert x - y  \vert^{2 - d} V(y) \omega(y) dy,
$$
we have that
\begin{equation}
\label{eq4.32}
(id + a\cF)(w^a) = -a\G
\end{equation}
where 
\begin{equation}
\nonumber
\G(x) = c_0 \int_{\mathbb{R}^d} |x-y|^{2-d}V(y)dy,
\end{equation}
 and $id$ means
the identity transformation. We wish to show the differetiability of $w^a$ with respect to $a>0$. 
This is clear heuristically because we have a candidate for $v^a:= \frac{\partial w^a}{\partial a}$;
 if we differentiate both sides of (4.4),
then  $v^a$ solves
\begin{equation}
\label{eq4.33}
(id + a\cF)(v^a) = -\G-\cF w^a=a^{-1}w^a.
\end{equation}
This provides us with a candidate for $\frac{\partial w^a}{\partial a}$, because the operator $id+a\cF$ has a bounded inverse
(see Section 6 of \cite{HR1}). The rigorous proof of the differentiability of $w^a$ goes as follows. First define $v^{a,h}=(w^{a+h}-w^a)/h$
and observe that $v^{a,h}$ satisfies
\begin{equation}
\label{eq4.34}
(id + a\cF)(v^{a,h}) = -\G-\cF w^{a+h}.
\end{equation}
We would like to show that $v^{a,h}$ has a limit in $\cH$, as $h\to 0$. One can readily show that the right-hand side of (4.6)
is bounded in $\cH$ because $|w^a(x)|\le C_2 a\min\{|x|^{2-d},1\}$ by the first part of the lemma. 
Hence $v^{a,h}$ stays bounded as $h\to 0$. If $v^a$ is any weak limit, then $v^a$ must satisfy (4.5).
Since (4.5) has a unique solution, the weak limit of $v^{a,h}$ exists. In \cite{HR1}, it is shown that $\cF$ is a compact operator. 
From this and (4.6), we can readily deduce that the strong limit of $v^{a,h}$ exists. As a consequence, $w^a$ is weakly differentiable in
$a$ and its derivative satisfies (4.5). Using Sobolev's inequalities and the fact that $V$ is H\"older continuous, we can deduce
by standard arguments  that indeed
$v^a$ is $C^2$ and satisfies
\begin{equation}\label{trieqn}
\triangle v^a=av^aV+(1+w^a)V.
\end{equation}
This means that $w^a(x)$ is continuously differentiable with respect to $(x,a)$.  

We now want to use (\ref{trieqn}) or equivalently (\ref{eq4.33})
to conclude that $a^{-1}w^a\le v^a\le 0$. In fact, by (\ref{eq4.33}), we have that $v^a=-a\cF v^a-a^{-1}w^a$, which implies that 
\[
|v^a(x)|\le c'_a c_0\int|x-y|^{2-d}dy +a^{-1}|w^a(x)|,
\]
where $c'_a$ is an upper bound for $|v^a(x)|$ with $x$ in the support of the function $V$. From this,
it is not hard to deduce that there exists a constant $c^{''}_a$ such that
\begin{equation}\label{vaineq}
|v^a(x)|\le c^{''}_a\max\{|x|^{2-d},1\}.
\end{equation}
In a similar fashion, we can show that there exists a constant $c^{'''}_a$ such that
\begin{equation}\label{delvaineq}
|\nabla v^a(x)|\le c^{'''}_a\max\{|x|^{1-d},1\}.
\end{equation}
 
We now demonstrate that $v^a\le 0$. Take a smooth function
$\varphi_{\d}: {\bR}
\to [0,\i)$ such that $ \varphi'_{\d}, \varphi_{\d} \ge 0$ and
\[
\varphi_{\d}(r) = \begin{cases}
0 &r \le 0, \\
 r &r \ge \d.
\end{cases}
\]
We then have
\begin{equation}
\label{eq6.18}
-\int_{\mathbb{R}^d} \varphi'_{\d}(v^a)|\nabla v^a|^2dx = 
\int_{\mathbb{R}^d} \varphi_{\d}(v^a)\D v^a dx =  \int_{\mathbb{R}^d}
V(1+w^a+av^a)\varphi_{\d}(v^a)dx,
\end{equation}
the second equality by (\ref{trieqn}). 
Integration by parts was performed in the first inequality:
we write the analogue of (\ref{eq6.18}) which is integrated over a bounded
set $\{x:|x|\le R\}$. We may obtain (\ref{eq6.18}) by sending $R\to\i$ but for this we need to make sure that the boundary contribution 
coming from the set $\{x:|x|=R\}$ goes away as $R\to\i$. This is readily achieved with the aid of (\ref{delvaineq}).
 Since $1+w^a\ge 0$ by the first part of the lemma, 
and $v^a \varphi_{\delta}(v^a) \geq 0$,
we deduce that the right-hand side of
\eqref{eq6.18} is non-negative.
Since the left-hand side is non-positive, we deduce that
\[
\int_{\mathbb{R}^d} \varphi'_{\d}(v^a)|\nabla v^a|^2dx = 
\int_{\mathbb{R}^d} V(1+w^a+av^a)\varphi_{\d}(v^a)dx = 0.
\]
We now send $\d \to 0$ to deduce
\[
0 = \int_{\mathbb{R}^d} |\nabla v^a|^21\!\!1(v^a \ge 0)dx = 
\int_{\mathbb{R}^d} V(1+w^a+av^a)v^a1\!\!1(v^a \ge 0)dx.
\]
As a result, on the set $A = \{x: v^a>0\}$ we have $\nabla v^a =
0$.  Hence
$v^a$ is constant on
each component $B$ of $A$.  But this constant can only be $0$
because on the boundary of $A$ we
have $v^a=0$.  This is impossible unless $A$ is empty. Hence,
$v^a \le 0$ everywhere.

It remains to prove that $v^a\ge a^{-1}w^a$. For this observe that if $\gamma^a=a^{-1}w^a-v^a$, then
\[
\triangle \g^a=aV\g^a+V(-w^a).
\]
We can now repeat the proof of $v^a\le 0$ to deduce that $\g^a\le 0$ because $-w^a\ge 0$. 
This completes the proof of the third part of the lemma.
\qed 

\noindent{\bf Proof of (2.12).}
 Note that
\begin{equation} \nonumber
\int^{T}_{0} \mathbb{E}_{N}\big|G_{z}^1-G_{0}^1 \big|(\mathbf q(t),t) dt \leq
\sum^{8}_{i=1} D_{i},
\end{equation}
where the first four of the $D_{i}$ are given by
\begin{eqnarray} \nonumber
D_{1} & = &\mathbb{E}_{N} \int^{T}_{0}dt \sum_{k, \ell \in I_{{\mathbf q}}}
\alpha(m_{k},m_{\ell})V_{\epsilon}(x_{k}-x_{\ell})\frac{m_{k}}{m_{k}+m_
{\ell}}\epsilon ^{2(d-2)} 
 \\ \nonumber
 && \sum_{i \in I_{{\mathbf q}}} \big| u^{\epsilon}(x_{k}-x_{i}+z;m_k+m_\ell,m_i) - u^{\epsilon}(x_{k}-x_{i};m_k+m_\ell,m_i)
 \big|\  \big| \hat J(x_{k},m_{k}+m_{\ell},x_{i},m_{i},t)
 \big| , \\ \nonumber
D_{2} & = &\mathbb{E}_{N} \int^{T}_{0}dt \sum_{k,\ell \in I_{{\mathbf q}}}
\alpha(m_{k},m_{\ell})V_{\epsilon}(x_{k}-x_{\ell})\frac{m_{\ell}}{m_{k}+m_
{\ell}}\epsilon^{2(d-2)} 
  \\ \nonumber
 && \sum_{i \in I_{{\mathbf q}}} \big| 
 u^{\epsilon}(x_{\ell}-x_{i}+z;m_k+m_\ell,m_i) - u^{\epsilon}(x_{\ell}-x_{i};m_k+m_\ell,m_i)
\big|\ \big| \hat J(x_{\ell},m_{k}+m_{\ell},x_{i},m_{i},t)
 \big| , \\ \nonumber
 D_{3} & = &\mathbb{E}_{N} \int^{T}_{0}dt \sum_{k,\ell \in I_{{\mathbf q}}}
\alpha(m_{k},m_{\ell})V_{\epsilon}(x_{k}-x_{\ell})\epsilon ^{2(d-2)}
 \\ \nonumber
 && \sum_{i \in I_{{\mathbf q}}} \big|  u^{\epsilon}(x_{k}-x_{i}+z;m_k,m_i  ) - u^{\epsilon}(x_{k}-x_{i};m_k,m_i  )
 \big| \ \big| \hat J(x_{k},m_{k},x_{i},m_{i},t)
 \big|  , \\ \nonumber
 \end{eqnarray}
 and
 \begin{eqnarray} \nonumber
  D_{4} & = &\mathbb{E}_{N} \int^{T}_{0}dt \sum_{k,\ell \in I_{{\mathbf q}}}
\alpha(m_{k},m_{\ell})V_{\epsilon}(x_{k}-x_{\ell})\epsilon ^{2(d-2)} 
 \\ \nonumber
 && \sum_{i \in I_{{\mathbf q}}} \big| u^{\epsilon}(x_{\ell}-x_{i}+z;m_\ell,m_i) - u^{\epsilon}(x_{\ell}-x_{i};m_\ell,m_i)
 \big|\ \big| \hat J(x_{\ell},m_{\ell},x_{i},m_{i})
 \big|  . \\ \nonumber
 \end{eqnarray}
 The other four terms each take the form of one of the above terms,
 the particles indices that appear in the arguments of the
 functions $u^{\epsilon}$ and $\hat J$ being switched, along with the mass pair
 labels for these functions.

 The estimates involved for each of the eight cases are in essence
 identical. We will examine the case of $D_{3}$. We write
 $D_{3}=D^{1}+D^{2}$, decomposing the inner $i$-indexed
 sum according to the respective index sets
 \begin{equation} \nonumber
 \big\{i \in I_{{\mathbf q}},i\neq k,\ell,|x_{k}-x_{i}| > \rho \big\}\,
 \text{and} \,  \big\{i \in I_{{\mathbf q}},i\neq k,\ell,|x_{k}-x_{i}| \leq \rho \big\}
 \end{equation}
 Here, $\rho$ is a positive parameter that satisfies the bound $\rho \geq 
\max \big\{2|z|+C_{0}
 \epsilon,2C_{0}\epsilon \big\}$. By the second part of Lemma 4.2, we
 have that
 \begin{equation} \nonumber
 D^{1} \leq \frac{c_0|z|\epsilon^{d-2}}{\rho^{d-1}}
 \mathbb{E}_{N}\int^{T}_{0}dt \sum_{k,\ell \in I_{{\mathbf q}}}
 \alpha(m_{k},m_{\ell})V_{\epsilon}(x_{k}-x_{\ell}),
 \end{equation}
 where we have also used the fact that the test function $\hat J$ is of compact support, and the fact that the total number of
 particles living at any given time is bounded above by $Z \epsilon^{2-d}$. 
From the bound on the collision that is provided by Lemma
 4.1, follows
 \begin{equation} \nonumber
 D^{1} \leq \frac{c_1|z|}{\rho^{d-1}}.
 \end{equation}

To bound the term $D^{2}$, note that by Lemma 4.2, the term $D^{2}$ is bounded above by
\begin{eqnarray} \nonumber
&& \mathbb{E}_{N} \int^{T}_{0} \epsilon
^{2(d-2)} \sum_{k,\ell \in I_{{\mathbf q}}}
\alpha(m_{k},m_{\ell})V_{\epsilon}(x_{k}-x_{\ell}) \\ \nonumber && \cdot \sum_{i \in I_{{\mathbf q}}} 1 \!\! 1
\big\{|x_{i}-x_{k}| \leq \rho \big\} \Big[
\big|u^{\epsilon}(x_{k}-x_{i}+z;m_k,m_i)\big|+\big|u^{\epsilon}(x_{k}-x_{i};m_k,m_i)
\big|\Big] \big| \hat{J}(x_{i},m_{i},x_k,m_k,t) \big|dt \\ \nonumber
&& \le c_1\mathbb{E}_{N} \int^{T}_{0} \epsilon
^{3(d-2)} \sum_{k,\ell \in I_{{\mathbf q}}}
\alpha(m_{k},m_{\ell})V^{\epsilon}(x_{k}-x_{\ell}) \\ \nonumber && \cdot \sum_{i \in I_{{\mathbf q}}} 1 \!\! 1
\big\{|x_{i}-x_{k}| \leq \rho ,\  \max \left\{ m_k,m_i,|x_k|,|x_i| \big\} \le L, \   m_k+m_i  \ge L^{-1}\right\} \nonumber \\
 & & \qquad \qquad \qquad \qquad \qquad \a'(m_k,m_i)\Big[
\big|x_k-x_{i}+z\big|^{2-d}+\big|x_{k}-x_{i}\big|^{2-d}\Big]dt , \nonumber
\end{eqnarray}
where $V^\e=\e^{2-d}V_\e$ and $L$ is chosen so that $\hat J(x,m,y,n)=0$ if any of the conditions 
$$
m+n\ge L^{-1},\ \ \  \max(m,n)\le L,\ \ \ \max(|x|,|y|)\le L,
$$
 does not hold. We note that if $m+k+m_i\ge L^{-1}$, then $\a'(m_k,m_i)\le c_2\a(m_k,m_i)$, for a constant $c_2$ that depends on $L$.
On the other hand, the conditions
$$
m_k\le L,\ \ \ m_i\le L,\ \ \ m_k {\text{ or }}m_i\ge \frac 12L^{-1},
$$ 
imply that for a constant $c_3=c_3(L)$,
$$
\a(m_k,m_{\ell})\a(m_k,m_i)\le c_3 \g_3(m_i)\g_3(m_\ell)\g_3(m_k),
$$
where we have used second part of Hypothesis 1.1.
We are now in a position to apply Corollary 3.1. 
For this we choose $k=3$ and 
\[
K(x_1,x_2,x_3)=V^\e(x_1-x_2)1 \!\! 1
\big\{ |x_{2}-x_{3}| \leq \rho,\ \  |x_{2}|,|x_{3}|\le L \big\}\Big[
\big|x_2-x_{3}+z\big|^{2-d}+\big|x_{2}-x_{3}\big|^{2-d}\Big] .
\]
As a result, $D^2\le D(z)+D(0)$ where $D(z)$ is given by
\begin{eqnarray} \nonumber
&&c_4\int V^\e(x_1-x_2)1 \!\! 1
\big\{|x_{2}-x_{3}| \leq \rho ,\ \ |x_{2}|,|x_{3}|\le L \big\}
\big|x_2-x_{3}+z\big|^{2-d} 
\prod_1^3 \left(\bar h_3*\l_3\right)(x_r)dx_r \\ \nonumber
&&\le c_5\int V^\e(x_1-x_2)1 \!\! 1
\big\{|x_{2}-x_{3}| \leq \rho ,\ \ |x_{2}|,|x_{3}|\le L\big\}
\big|x_2-x_{3}+z\big|^{2-d} dx_1dx_2dx_3\\ \nonumber
&&\le c_6\int_{|a|\le \rho}|a+z|^{2-d}da\le c_7(\rho+|z|)^2,
\end{eqnarray}
where, for the first inequality, we used  Hypothesis 1.2(ii). 
Combining these estimates yields
\begin{equation} \nonumber
D_{3} = D^{1}+ D^{2} \leq c_1\frac{|z|}{\rho^{d-1}}+ c_7\big(\rho +
|z| \big)^{2}.
\end{equation}
Making the choice $\rho=|z|^{\frac{1}{d+1}}$ leads to the inequality
$D_{3} \leq c_8 |z|^{\frac{2}{d+1}}$. Since each of the cases of
$\Big\{D_{i}:i \in \{1, \ldots ,8\} \Big\}$ may be treated by a
nearly verbatim proof, we are done.
\qed

\textbf {Proof of (2.15)}.
Setting $\mathbb{L} = \mathbb{A}_0 + \mathbb{A}_c$,
the process
\begin{equation} \nonumber
M_{z}(T)=X_{z}({\mathbf q}(T),T)-X_{z}({\mathbf q}(0),0)-  \int_{0}^{T} \Big(
\frac{\partial}{\partial t}+\mathbb{L} \Big)X_{z}({\mathbf q}(t),t) dt
\end{equation}
is a martingale which satisfies
\begin{equation} \nonumber
\mathbb{E}_{N}\big[M_{z}(T)^{2}\big] =
\mathbb{E}_{N}\int_{0}^{T} \left(
\mathbb{L}X_{z}^2- 2X_{z} \mathbb{L}X_{z} \right)({\mathbf q}(t),t) dt
=\sum_{i=1}^{3}
\mathbb{E}_{N}~\int_{0}^{T} A_{i}({\mathbf q}(t),t) dt,
\end{equation}
where
$$
A_{1}({\mathbf q},t) =  2 \epsilon^{4(d-2)} \sum_{i \in I_{{\mathbf q}}} d(m_{i}) 
\Big[\nabla _{x_{i}} \sum_{j \in I_{{\mathbf q}}} u^{\epsilon} (x_{i}-x_{j}+z;m_i,m_j)
\hat J(x_{i},m_{i},x_{j},m_{j},t) \Big]^{2},
$$
and
$$
A_{2}({\mathbf q},t) = 2 \epsilon^{4(d-2)} \sum_{j \in I_{{\mathbf q}}} d(m_{j}) 
\Big[\nabla _{x_{j}} \sum_{i \in I_{{\mathbf q}}} u^{\epsilon} (x_{i}-x_{j}+z;m_i,m_j)
\hat J(x_{i},m_{i},x_{j},m_{j},t) \Big]^{2},
$$
while $A_{3}({\mathbf q},t)$ is given by
\begin{eqnarray} \label{eq4.3a}
& &\epsilon^{4(d-2)}  \sum_{i,j \in I_{{\mathbf q}}} \alpha(m_{i},m_{j})
\epsilon^{-2} V_\e( {x_{i}-x_{j}}) \\
\nonumber && 
\quad\Big\{ \sum_{k \in I_{{\mathbf q}}} \Big[
\frac{m_{i}}{m_{i}+m_{j}} u^{\epsilon}(x_{i}-x_{k}+z;m_i+m_j,m_k)
\hat J(x_{i},m_{i}+m_{j}, x_{k},m_{k},t) \\ 
\nonumber && \quad\quad+
\frac{m_{i}}{m_{i}+m_{j}} u^{\epsilon}(x_{k}-x_{i}+z;m_k,m_i+m_j) \hat J(x_{k},m_{k},x_{i},m_{i}+m_{j},t) \\ \nonumber
&&\quad\quad+\frac{m_{j}}{m_{i}+m_{j}} u^{\epsilon}(x_{j}-x_{k}+z;m_i+m_j,m_k)
\hat J(x_{j},m_{i}+m_{j}, x_{k},m_{k},t) \\ 
\nonumber
&&\quad\quad
+\frac{m_{j}}{m_{i}+m_{j}} u^{\epsilon}(x_{k}-x_{j}+z;m_k,m_i+m_j)
\hat J(x_{k},m_{k}, x_{j},m_{i}+m_{j},t) \\ 
\nonumber
&&\quad\quad-u^{\epsilon}(x_{i}-x_{k}+z;m_i,m_k) \hat J(x_{i},m_{i},
x_{k},m_{k},t) \\ 
\nonumber
&&\quad\quad-u^{\epsilon}(x_{k}-x_{i}+z;m_k,m_i) \hat J(x_{k},m_{k},
x_{i},m_{i},t) \\ 
\nonumber
&&\quad\quad-u^{\epsilon}(x_{j}-x_{k}+z;m_j,m_k) \hat J(x_{j},m_{j},
x_{k},m_{k},t) \\ 
\nonumber
&&\quad\quad-u^{\epsilon}(x_{k}-x_{j}+z;m_k,m_i) \hat J(x_{k},m_{k},
x_{j},m_{j},t)  \Big] \\ 
\nonumber
&&\quad\quad-u^{\epsilon}(x_{i}-x_{j}+z;m_i,m_j) \hat J(x_{i},m_{i},
x_{j},m_{j},t) \Big\} ^{2} 
\end{eqnarray}
We now bound the three terms. Of the first two, we treat only
$A_{1}$, the other being bounded by an identical argument. By
multiplying out the brackets appearing in the definition of $A_{1}$,
and using $\sup_{m \in (0,\infty)}d(m) < \infty$, (which is assumed by Hypothesis $1.1$),
we obtain that $A_{1}\le A_{11}+A_{12}$ with
\begin{eqnarray} \nonumber
A_{11}& =&c_0 \epsilon^{4(d-2)} \sum_{i,j,k \in I_{{\mathbf q}}} 
\left|u_x^\e\left({x_{i}-x_{j}+z} ;m_i,m_j \right)\right |
\left| u_x^\e\left({x_{i}-x_{k}+z} ;m_i,m_k \right)\right |\\ \nonumber
&& \quad\quad \quad\quad\quad\cdot|\hat J(x_{i},m_{i},x_j,m_j,t)|| \hat J(x_{i},m_{i},x_k,m_k,t) |
\\
\nonumber  A_{12}& =&
c_0 \epsilon^{4(d-2)}  \sum_{i,j,k \in I_{{\mathbf q}}} 
\left| u^\e \left({x_{i}-x_{j}+z} ;m_i,m_j \right)\right| 
\left| u^\e\left({x_{i}-x_{k}+z} ;m_i,m_k \right)\right|\\
\nonumber
&&\quad\quad \quad\quad \quad \cdot|\hat J_x(x_{i},m_{i},x_{j},m_{j},t)|  |\hat J_x(x_{i},m_{i},x_{k},m_{k},t)|.    
\end{eqnarray}
Let us assume that $z=0$ because this will not affect our arguments. We bound the term $A_{11}$ with the aid of 
Corollary 3.1 and Lemma 4.2. The term $A_{12}$ can be treated likewise.  To bound $A_{11}$, first observe even though
 $i$ and $j$ are distinct,
$k$ and $j$ can coincide. Because of this, let us write $A_{11}=A_{111}+A_{112}$ where $A_{111}$ represents the case of distinct $i,j$
and $k$. We only show how to bound $A_{111}$ where the correlation bound in the case of $k=3$ is used. The term
$A_{112}$ can be treated in the similar fashion with the aid of Corollary 3.1 when $k=2$. Since 
$\hat J(x,m,y,n)\neq 0$ implies that $m,n,|x|,|y|\le L$ and $m+n\ge L^{-1}$. 
Using second part of Hypothesis 1.1,  we can find a constant $c_1=c_1(L)$ such that
$$
\a(m_i,m_j)\a(m_i,m_k)\le c_2\g_3(m_i)\g_3(m_j)\g_3(m_k),
$$ 
whenever
$$
m_i,m_j,m_k\le L,\ \ \ m_i+m_j,m_i+m_k\ge L^{-1}.
$$
As a result, we may apply Corollary 3.1 with $k=3$ and 
\[
K(x_1,x_2,x_3)= \e^{d-2}|x_1-x_2|^{1-d}|x_1-x_3|^{1-d}1\!\!1(|x_1|,|x_2|,|x_3|\le L),
\]
to deduce
\[
A_{111}  \le c_2 \e^{d-2} \int |x_1-x_2|^{1-d}|x_1-x_3|^{1-d}1\!\!1(|x_1|,|x_2|,|x_3|\le L)
\prod_{r=1}^3 \left(\bar h_3*\l_3\right)(x_r)dx_r .
\]
Note that $K$ is an unbounded function and Corollary 3.1 can not be applied directly. However we can approximate $K$ with a 
sequence of bounded functions and pass to the limit. From this and  Hypothesis 1.2, we deduce
\[
A_{11}\le c_3 \e^{d-2} \int |x_1-x_2|^{1-d}|x_1-x_3|^{1-d}1\!\!1(|x_1|,|x_2|,|x_3|\le L)dx_1dx_2dx_3=c_4 \e^{d-2} .
\]
This and an analogous argument that treats the terms $A_{112}$, $A_{12}$ and $A_2$ lead to the conclusion that
\begin{equation}
A_1+A_2\le c_4 \e^{d-2} .
\end{equation}

We must treat the third term, $A_{3}$. An application of the
inequality
\begin{equation} \nonumber
(a_{1}+ \ldots +a_{n})^{2} \leq n(a_{1}^{2}+ \ldots +a_{n}^{2})
\end{equation}
to $A_{3}$, given in \eqref{eq4.3a}, implies that
\begin{equation} \label{eq4.13}
A_{3}(\mathbf q,t) \leq 9 \epsilon^{4(d-2)} \sum_{i,j \in I_{\mathbf q}} \alpha
(m_{i},m_{j}) V_\e({x_{i}-x_{j}})\Big[ \sum^{8}_{n=1}\big(\sum_{k \in I_{\mathbf q} }Y_{n} \big)^{2}
+Y_{9}^{2} \Big]=:\sum_{i=1}^9A_{3i},
\end{equation}
where $Y_{1}$ is given by
\begin{equation} \nonumber
\frac{m_{i}}{m_{i}+m_{j}} u^{\epsilon} (x_{i}-x_{k}+z;m_i+m_j,m_k)
\hat{J}(x_{i},m_{i}+m_{j},x_{k},m_{k},t),
\end{equation}
and where $\{ Y_{i}: i \in \{ 2, \ldots ,8\} \}$ denote the other
seven expressions in \eqref{eq4.3a} that appear in a sum over $k \in
I_{q}$, while $Y_{9}$ denotes the last term in \eqref{eq4.3a} that
does not appear in this sum. There are nine cases to consider. The
first eight are practically identical, and we treat only the fifth.
Let us again assume that $z=0$ because this will not affect our arguments. Note that
\begin{eqnarray} \nonumber
 A_{35}&=&\epsilon ^{4(d-2)} \sum_{i,j \in I_{\mathbf q}} \alpha (m_{i},m_{j})
V_{\epsilon} (x_{i}-x_{j}) \big(\sum_{k \in I_{\mathbf q}} Y_{5} \big)^{2}
\\ \nonumber
& =&  \epsilon ^{5(d-2) } \sum_{i,j \in I_{\mathbf q}}\alpha (m_{i},m_{j})
V^{\epsilon} (x_{i}-x_{j})   \\  && \ \Bigg[ \sum_{k,l \in I_{\mathbf q}} u^\e
\left( {x_{i}-x_{k}};m_i,m_k \right) u^\e \left(
{x_{i}-x_{l}} ;m_i,m_l\right)\hat J(x_{i},m_{i},
x_{k},m_{k},t) \hat{J}(x_i,m_i,x_{l},m_{l},t) \Bigg].
\nonumber
\end{eqnarray}
In the sum with indices involving $k,l \in I_{\mathbf q}$, we permit the
possibility that these two may be equal, though they must be
distinct from each of $i$ and $j$ (which of course must themselves
be distinct by the overall convention). 
Let us write $A_{35}=A_{351}+A_{352}$, where 
$A_{351}$ corresponds to the case when all the indices $i,j,k$ and $l$ are distinct and $A_{352}$ corresponds to the remaining cases.
Again, our assumption on $\a$ as in Hypothesis 1.2 would allow us 
to treat the term $A_{351}$ with the aid of Corollary 3.1.  This time $k=4$ and  our bound on $u$
 given in the first part of Lemma 4.2 suggests the following choice for $K$:
\[
K(x_1,\dots,x_4)= \e^{d-2}V^\e(x_1-x_2)|x_1-x_3|^{2-d}|x_1-x_4|^{2-d}1\!\!1(|x_1|,|x_2|,|x_3|,|x_4|\le L).
\]
Note that $K$ is an unbounded function and Corollary 3.1 can not be applied directly. However we can approximate $K$ with a 
sequence of bounded functions and pass to the limit. From Corollary 3.1 and  Hypothesis 1.1 on the initial data 
we deduce that the expression $\int_0^T A_{351}dt$ is bounded above by 
\[
c_5\e^{d-2}\int V^\e(x_1-x_2)|x_1-x_3|^{2-d}|x_1-x_4|^{2-d}1\!\!1(|x_1|,|x_2|,|x_3|,|x_4|\le L)dx_1\dots dx_4
=c_6\e^{d-2}.
\]
A similar reasoning applies to $A_{352}$, except that Corollary 3.1 in the case of $k=3$ would be employed. Hence,
\begin{equation}
\sum_{i=1}^8A_{3i}\le c_7\e^{d-2}.
\end{equation}

We now treat the ninth term, as they are classified in
\eqref{eq4.13}. It takes the form
\[
 \epsilon^{4d-8 }  \sum_{i,j \in I_{\mathbf q}} \alpha(m_{i},m_{j})
V_{\epsilon}(x_{i}-x_{j}) 
  u^{\epsilon} (x_{i}-x_{j}+z;m_i,m_j)^{2} \hat J(x_{i},m_{i},x_{j},m_{j},t)^{2}.
\]
This is bounded above by
\begin{equation} \nonumber
c_8 \epsilon^{2d-4} \sum_{i,j \in I_{\mathbf q}} \alpha(m_{i},m_{j})
V_\e(x_{i}-x_{j}),
\end{equation}
because $u^{\epsilon} \leq c_9 \epsilon^{2-d}$ by the first part of
Lemma 4.2. The expected value of the integral on the interval of
time $[0,T]$ of this last expression is bounded above by
\begin{equation} \nonumber
c_7\epsilon^{2d-4} \mathbb{E}_{N} \int_{0}^{T}  \sum_{i,j \in I_{q}}
\alpha(m_{i},m_{j}) V_{\epsilon}(x_{i}-x_{j})dt \leq c_{10} \epsilon
^{d-2}.
\end{equation}
where we used Lemma 4.1 for the last inequality. This, (4.12), (4.13)  and (4.14) complete the
proof of (2.15).
\qed

\end{section}
\setcounter{equation}{0}

\begin{section}{ Bounds on the Macroscopic Densities}
In this section we show how Corollary 3.1 can be used to obtain certain bounds on the macroscopic densities. 
These bounds will be used for the derivation of the macroscopic equation. Recall that
\begin{equation} \nonumber
g^{\epsilon} (dx,dn,t) = \epsilon^{d-2} \sum_{i}
\delta_{(x_{i}(t),m_{i}(t))}(dx,dn),
\end{equation}
and that the law of
\begin{equation}\nonumber
\mathbf{q} \mapsto g^{\epsilon} (dx,dn,t)
\end{equation}
induces a probability measure $\mathcal{P}^{\epsilon}$ on the
 space $\cX$. 
The main result is Theorem 5.1.

\begin{theorem} Let $\cP$ be a limit point
of $\cP^\e$. The following statements are true:
\begin{itemize}
\item
1. For every positive $L_1$, and $k\in\{2,3,4\}$,
\begin{equation} \label{eq5.31}
\sup_\d \int_{\cX}\int_0^\i \int_{|x|\le L_1} \Big[\int_{0}^\i\int \xi^{\delta} (x-y) \g_k(n)
g(dy,dn,t) \Big]^{k}dxdt d\mathcal{P} <\i,
\end{equation}
where $\xi^{\delta}(x)=\delta^{-d}\xi\Big(\frac{x}{\delta} \Big)$, with $\xi$  a  nonnegative smooth function of compact support
satisfying $\int\xi=1$. 
\item 
2. We have $g(dx,dn,t)=f(x,t,dn)dx$ for almost all $g$ with respect to the probability measure $\cP$.
\item
3. For every continuous $R$ of compact support and positive $L$,
\begin{eqnarray} \label{eq5.4}
&& \lim_{\delta \rightarrow 0} \int\Big|\int_0^T \int ^{L}_{L^{-1}} \int
^{L}_{L^{-1}} \int   R(x,m,n,t)    f^{\delta}(x,t,dm)
f^{\delta}(x,t,dn) dx dt\\ \nonumber & & -\int_0^T \int
^{L}_{L^{-1}}\int ^{L}_{L^{-1}} \int R(x,m,n,t) f(x,t,dm)
f(x,t,dn) dxdt\Big|d\cP=0,
\end{eqnarray}
 where
\begin{equation}\label{fdeqn}
f^{\delta}(x,t,dn)=\int\xi^{\delta} (x-y)g(dy,dn,t).
\end{equation}
\end{itemize}
\end{theorem}

{\textbf {Proof.}}
Fix $x\in \bR^d$ and
choose
\begin{equation} \nonumber
K(y_{1}, \ldots,y_{k})= \prod ^{k}_{r=1} \xi^{\delta}(x-y_{r}),
\end{equation}
in Corollary 3.1. The right-hand side of (3.7) equals
\begin{equation} \nonumber
\int 
\prod^{k}_{r=1}  \xi^{\delta}(x-x_{{r}})\bar h_k*\l_k(x_r)dx_r,
\end{equation}
which, by the second part of Hypothesis 1.2, is bounded by a constant $c_1(L_1)$ when $k=2,3,4$,   and 
$|x|\le L_1$. As a result,
\begin{equation} \label{eq5.2}
\mathbb{E}_N\int_0^\i\int_{|x|\le L_1} \epsilon^{k(d-2)} \sum_{i_{1}, \ldots , i_{k}} \prod ^{k}_{r=1}
\xi^{\delta} (x-x_{i_{r}}(t)) \g_k(m_{i_{r}}(t))  dxdt\\
 \leq  c_1(L_1) 
\end{equation}
for a constant $c_1(L)$ which is independent of $\delta$ and $\e$. Here we are assuming that the indices
$i_1,\dots,i_k$ are distinct. Note that if we allow non-distinct 
indices in the summation, then the difference would go to $0$ as $\e\to 0$ because the summation is multiplied
by $\e^{k(d-2)}$ while the number of additional terms is of order $O(\e^{(k-1)(2-d)})$. As a consequence, we can
use (5.4) to deduce (5.1).

Recall that the function $\g_k$ is a positive continuous function. From this and (5.1), one can readily deduce  part 2.

It remanis to establish part 3.  First observe that by (5.1) and the posivity of $\g_4$,
\begin{equation} \label{eq5.311}
\sup_{\delta} \int\int_0^T \int_{|x|\le L_1} \left[\int_{L^{-1}}^L f^\d(x,t,dn)\right]^{4}dxdt \mathcal{P}(dg) \leq c_2(L_1,L).
\end{equation}
Because of this, 
it suffices to prove that
\begin{eqnarray} \nonumber
&& \lim_{\delta \rightarrow 0} \int_0^T\int \int^{L}_{L^{-1}}
\int^{L}_{L^{-1}} {R}_{p}(x,m,n,t)  f^{\delta}(x,t,dm)
f^{\delta}(x,t,dn) dx dt\\ \nonumber && ~= \int_0^T \int
^{L}_{L^{-1}}\int ^{L}_{L^{-1}} \int {R_p}(x,m,n,t)  f(x,t,dm)
f(x,t,dn)  dx.
\end{eqnarray}
for each $p$, provided that  $\lim_{p \rightarrow \infty}
{R}_{p}(x,m,n,t)={R}(x,m,n,t)$, uniformly for $m,n \in
[L^{-1},L]$, $|x|\le L_1$ and $t\le T$. By approximation, we may assume that $R$ is
of the form $R(x,m,n,t)= \sum^{\ell}_{i=1} J_{1}^\ell(x,t)J_{2}^\ell(m)J_3^\ell(n).$ Hence it suffices to
establish \eqref{eq5.4} for ${R}$ of the form $R(x,m,n,t)=  J_{1}(x,t)J_{2}(m)J_3(n).$ 
 But now the left-hand side of
\eqref{eq5.4} equals
\begin{equation} \nonumber
\lim_{\delta \rightarrow 0} \int_0^T \int \left[\int^{L}_{L^{-1}}
 J_{2}(m)f^{\delta}(x,t,dm) \right] \left[
\int^{L}_{L^{-1}}  J_3(n) f^{\delta}(x,t,dn) \right]
J_{1}(x,t)dxdt.
\end{equation}
We note that
\begin{equation} \nonumber
\int^{L}_{L^{-1}}  J_{2}(m)f^{\delta}(x,t,dm ) =
\left( \int^{L}_{L^{-1}} J_{2}(m)f(\cdot,t,dm)
 \right)\ast_{x} \xi ^{\delta}(x).
\end{equation}
converges almost everywhere to
\begin{equation} \nonumber
\int^{L}_{L^{-1}} J_{2}(m) f(x,t,dm) .
\end{equation}
The same comment applies to $\int^{L}_{L^{-1}} J_3(n) 
f^{\delta}_{n}(x,t) dn.$ From this and (5.5) we deduce (5.2).
\qed

 \setcounter{equation}{0}
\begin{section}{Deriving the PDE}
We wish to derive  (1.6) from the identity (2.1). There is a technical issue we need to settle
first: in (2.2), the function $\hat J(x,m,y,n,t)$ does not have a compact support with respect to $(m,n)$,
 even if $J$ is of compact support.
Recall that in Theorem 2.1 we have assumed that $\hat J$ is of compact support. Lemma 6.1 settles this issue.
\begin{lemma}
There exists a constant $C_4$ independent of $\e$ such that
\begin{equation}
\bE_N\int_0^T\epsilon^{2(d-2)} \sum_{i,j \in  I_{{\mathbf q}}}
\alpha(m_{i}(t),m_{j}(t)) V_{\epsilon} (x_{i}(t)-x_{j}(t)) m_i(t)m_j(t)dt\le C_4.
\end{equation}
Moreover,
\begin{equation}
\lim_{L\to \i}\sup_\e\bE_N\int_0^T\epsilon^{2(d-2)} \sum_{i,j \in  I_{{\mathbf q}}}
\alpha(m_{i}(t),m_{j}(t)) V_{\epsilon} (x_{i}(t)-x_{j}(t))1\!\!1( \min \big\{ m_i(t) , m_j(t) \big\} \le L^{-1} )dt=0.
\end{equation}
\end{lemma}

\textbf {Proof.} Let us take a smooth function $J: \bR^d 
\to [0,\i)$ and set
\begin{equation}
\label{eq3.2a}
H(x) = c_0(d) \int \frac {J(y)}{|x-y|^{d-2}} dy
\end{equation}
with $c_0(d) = (d-2)^{-1}\omega_d^{-1}$ with $\omega_d$ denoting the surface 
area of the unit sphere in $\bR^d$.  Note that $H \ge 0$ and $-\D H = 
J$.  Let $\psi:(0,\i)\times (0,\i)\to [0,\i)$ be a continous symmetric function and set
\begin{equation}
\label{eq3.3}
X_N(\mathbf q) = {\e}^{2(d-2)} \sum_{i,j \in I_{\mathbf q}} H(x_i-x_j)\psi(m_i,m_j)
\end{equation}
We have
\begin{eqnarray}
\label{eq6.4}
-\bE_N \int_0^T \cA_cX_N({\mathbf q}(s))ds - \bE_N \int_0^T \cA_0X_N({\mathbf q}(s))ds &= &
\bE_N X_N({\mathbf q}(0)) - \bE_N X_N({\mathbf q}(T))\nonumber\\
& \le& \bE_N X_N({\mathbf q}(0)),
\end{eqnarray}
where
\[
\cA_0X_N({\mathbf q})= -\e^{2(d-2)} \sum_{i,j \in I_{\mathbf q}} J(x_i-x_j)\psi(m_i,m_j)(d(m_i)+d(m_j)),
\]
and $\cA_cX_N({\mathbf q})=Y_1(\mathbf q)+Y_2(\mathbf q)$, with 
\begin{eqnarray*}
Y_1(\mathbf q)&= & -{\e}^{2(d-2)}  \sum_{i,j \in I_{\mathbf q}} \a(m_i,m_j) 
V_{\e}(x_i-x_j) \psi(m_i,m_j)H(x_i-x_j) \\
Y_2(\mathbf q)&= & {\e}^{2(d-2)}  \sum_{i,j,k \in I_{\mathbf q}} \a(m_i,m_j) 
V_{\e}(x_i-x_j)\G(x_i,x_j,x_k,m_i,m_j,m_k),
\end{eqnarray*}
where
\begin{eqnarray*}
\G(x_i,x_j,x_k,m_i,m_j,m_k)&=& \left[ \frac {m_i}{m_i+m_j} \psi(m_i+m_j,m_k)-\psi(m_i,m_k)\right]H(x_i-x_k)\\
&&+\left[ \frac{m_j}{m_i+m_j}  \psi(m_i+m_j,m_k)-\psi(m_j,m_k)\right]H(x_j-x_k)\\
&& +\left[ \frac {m_i}{m_i+m_j} \psi(m_k,m_i+m_j)-\psi(m_k,m_i)\right]H(x_k -x_i)\\
&&+\left[ \frac{m_j}{m_i+m_j}  \psi(m_k,m_i+m_j)-\psi(m_k,m_j)\right]H(x_k-x_j).
 \end{eqnarray*}
We consider two examples for $\psi$. As the first example, we choose $\psi(m,n)=mn$. This yields $Y_2 = 0$. 
We find that
\begin{equation}
\label{eq6.5}
\sup_N\bE_N \int_0^T Y_1({\mathbf q}(s))ds \le \bE_N X_N({\mathbf q}(0)).
\end{equation}
The hope is that a 
suitable choice of $J$ would yield the desired assertion 
(6.1).  For this, we simply choose $J(x) = {\e}^{-d}A\left( \frac {x}{\e} \right)$ where $A$ is a smooth non-negative 
function of compact support.  We then have that $H(x) = \e^{2-d}B 
\left( \frac {x}{\e} \right)$ where $\D B = -A$.  As a result,
\begin{equation}
\label{eq6.6}
Y_1({\mathbf q}) = {\e}^{d-2} \sum_{i,j \in I_{\mathbf q}} V_{\e}(x_i-x_j) B\left( 
\frac {x_i-x_j}{\e} \right) m_im_j\a(m_i,m_j)
\end{equation}
with
\begin{equation}\nonumber
B(x) = c_0(d) \int \frac {A(y)}{|x-y|^{d-2}} dy.
\end{equation}
Recall that the support of $V$ is contained in the set $y$ with $|y|\le C_0$.
If we choose $A$ so that
\[
1\!\!1(|y|\le 3C_0)\le A(y) \le  1\!\!1(|y|\le 4C_0),
\]
 then, for $|x| \le C_0$,
\begin{equation}
\label{eq6.8}
B(x) \ge c_0(d) \int_{3C_0\ge |y| \ge 2{C_0}} \frac {dy}{|x-y|^{d-2}} \le c_0(d) C_0^{2-d}
\int_{3C_0\ge |y| \ge 2{C_0}}  {dy} =: \tau_0 > 0.  \nonumber
\end{equation}
 On the other hand, if $|x|\le 5 C_0$, then
\begin{equation}
\label{eq6.8a}
B(x) \le c_0(d) \int_{|x-y|\le 9C_0} \frac {dy}{|x-y|^{d-2}} = \frac 12c_0(d) \o_d  (9C_0)^{2}.
\end{equation}
and if $|x|\ge 5C_0$, then
\[
B(x) \le c_0(d) \left|\frac {4x}{5}\right|^{2-d} \int_{C_0\ge |y|}  {dy} =c_1\left| x\right|^{2-d} .
\]
From this, (6.8) and the third part of Hypothesis 1.2, we learn that the right-hand side of (6.6)  is uniformly bounded in $\e$.
 This completes the proof of (6.1).

As for (6.2), we choose $\psi(m,n)=1\!\!1(m\le \d)+1\!\!1(n\le\d)$. This time we have that $Y_2\le 0$.
Such a function $\psi$ is not continuous. But by a simple approximation procedure we can readily see that (6.5) is valid 
for such a choice. By the third part of Hypothesis 1.2 on the initial data, we know that 
\[
\int_0^\i\int h_n(x)\hat h(y)|x-y|^{2-d}dxdydn<\i.
\]
From this we learn that
\[
\lim_{\d\to 0}\int_0^\d\int h_n(x)\hat h(y)|x-y|^{2-d}dxdydn=0,
\]
whence 
\[
\lim_{\d\to 0}\sup_N\bE_NX_N(\mathbf q(0))=0.
\]
This and (6.5) imply (6.2).
\qed

{\textbf {Proof of Theorem 1.1.}} \noindent{\textbf{Step 1.}} 
We take a  smooth test function $J$ of compact support in $\bR^d\times (0,\i)\times [0,\i)$
and study the decomposition (2.1). Firstly, we 
show that the martingle term goes to $0$. 
 The term $M_{T}$ is a martingale satisfying
\begin{equation} \nonumber
\mathbb{E}_{N}\big[ M_{T}^{2}\big]=\mathbb{E}_{N}\int_0^T\left(\mathbb{L} Y^2-2Y\mathbb{L} Y\right)(\mathbf q(t),t)dt
= \mathbb{E}_{N} \int ^{T}_{0}
A_{1}({\mathbf q}(t),t) dt + \mathbb{E}_{N} \int^{T}_{0}A_{2}({\mathbf q}(t),t) dt ,
\end{equation}
where $A_{1}({\mathbf q},t)$ and $A_{2}({\mathbf q},t)$ are respectively set equal to
\begin{equation} \nonumber
A_{1}({\mathbf q},t)=\epsilon^{2(d-2)} \sum_{i \in I_{{\mathbf q}}} d(m_{i}) |
 J_x (x_{i},m_{i},t)|^{2} ,
\end{equation}
and
\[
A_{2}({\mathbf q},t) =  \epsilon^{2(d-2)} \sum_{i \in  I_{{\mathbf q}}}
\alpha(m_{i},m_{j}) V_{\epsilon} (x_{i}-x_{j}) \hat J(x_i,m_i,x_j,m_j,t)^2.
\]
 We  can readily show
\begin{eqnarray} \label{eq6.2}
A_{1}({\mathbf q},t)& \leq &c_1\epsilon^{2(d-2)} \sum_{i \in I_{{\mathbf q}}} d(m_{i})
\leq c_2 \epsilon^{d-2}, \\ \label{eq6.3}
\mathbb{E}_{N} \int^{T}_{0} A_{2}({\mathbf q}(t),t)dt & \leq &
c_3 \mathbb{E}_{N} \int^{T}_{0}\epsilon^{2(d-2)} \sum_{i,j \in I_{{\mathbf q}}} \alpha(m_{i},m_{j})
V_{\epsilon}(x_{i}-x_{j})dt \leq c_4 \epsilon^{d-2},
\end{eqnarray}
where we have Lemma 4.1 in the last inequality. From these
inequalities, we deduce that the martingale
tends to zero, in the $\epsilon \downarrow 0$ limit.

\textbf{Step 2.} 
We rewrite the terms of (2.1) in terms of the empirical
measures. We have that
\begin{equation} \label{eq6.4a}
Y({{\mathbf q}}(t),t) = \int_0^\i \int_{\mathbb{R}^{d}} J(x,n,t)
g(dx,dn,t), 
\end{equation}
and that
\begin{equation} \label{eq6.5a}
\int^{T}_{0} \left(\frac{\partial}{\partial t}+ \mathbb{A}_{0}\right)Y({{\mathbf q}}(t),t)dt = \int^{T}_{0}
\int_{0}^{\infty}\int_{\mathbb{R}^{d}} \left(\frac{\partial}{\partial t}+ d(n) \triangle_{x}\right) J(x,n,t)
g(dx,dn,t) .
\end{equation}
Furthermore, by Theorem 2.1 and Lemma 6.1,
\begin{equation}
\int^{T}_{0}   \mathbb{A}_{c}Y({\mathbf q}(t),t)dt = \int_0^T \Gamma^\d_L(\mathbf q(t),t)dt + Err^1(\e,L)+Err^2(\e,\d,L),
\end{equation}
where $T$ is large enough so that $J(\cdot,\cdot,t)=0$ for $t\ge T$,  the expression $\Gamma^\d_L(\mathbf q,t)$ is given by
\[
\iint
\int_{L^{-1}}^L \int_{L^{-1}}^L \alpha(m,n)U_{n,m}^\epsilon(w_1-w_2)f^\delta(w_1,dm;\mathbf q)
f^\delta(w_2,dn;\mathbf q)\hat J(w_1,m,w_2,n,t)dw_1dw_2,
\]
and
$$
\lim_{L\to\i}\sup_\e\bE_N|Err^1(\e,L)|=0,\ \ \ \lim_{\d\to 0}\limsup_{\e\to 0}\bE_N|Err^2(\e,\d,L)|=0.
$$
We note that if  we replace $f^\delta(w_2,dn;\mathbf q)\hat J(w_1,m,w_2,n,t)$ with 
$f^\delta(w_1,dn;\mathbf q)\hat J(w_1,m,w_1,n,t)$,
then we produce an error which is of order $O(L\d^{-\d-1}\e)$, which goes to $0$ because we send $\e\to 0$ first.
As a result, (6.13) equals
\begin{eqnarray} \nonumber
&& \int^{\i}_{0}  \int_{\mathbb{R}^{d}}  \int_{L^{-1}}^{L} \int_{L^{-1}}^{L} \b(m,n) (g \ast_x
\xi^{\delta})(x,t,dm) (g \ast_x \xi^{\delta})(x,t,dn)  \tilde J(x,m,n,t)  dx dt\\
\nonumber && ~~~ + Err^1(\e,L)+  Err^3(\e,\d,L) ,
\end{eqnarray}
where 
\[
\lim_{\d\to 0}\limsup_{\e\to 0}\mathbb{E}_{N}| Err^3(\epsilon,\delta,L)| =0.
\]
By passing to the limit in low $\epsilon$, we find
that any weak limit $\mathcal{P}$ is concentrated on the space of
measures $g(dx,dn,t)dt$ such that,
\begin{eqnarray} \nonumber
& & \int_{0}^{\infty}\int_{\mathbb{R}^{d}}h_n(x)J(x,n,0)dxdn+\int_{0}^{\infty} \int^{\i}_{0} \int_{\mathbb{R}^{d}} 
g(dx,dn,t)\left(\frac {\partial}{\partial t}+ d(n)\triangle_x\right) J(x,n,t)dt
\\  && +\int^{\i}_{0}\int_{\mathbb{R}^{d}}   \int_{L^{-1}}^{L} \int_{L^{-1}}^{L}  \b(m,n) (g \ast_x
\xi^{\delta})(x,t,dm) (g \ast_x \xi^{\delta})(x,t,dn)  \tilde J(x,m,n,t)  dx dt\\
\nonumber && 
+Err^4(L)+Err^5(\d)=0,
\end{eqnarray}
where  the $\mathcal{P}$-expectation of
$|Err^5(\delta)|$ goes to zero as $\delta \downarrow 0$, and the $\mathcal{P}$-expectation of
$|Err^4(L)|$ goes to zero as $L\to \i$.  From Theorem 5.1 we know that $g(dx,dn,t)=f(x,t,dn)dx$, $\cP$-almost surely 
and that by (5.2) we can replace $g*_x\xi$ with $f$. Hence
\begin{eqnarray} \nonumber
& & \int_{0}^{\infty}\int_{\mathbb{R}^{d}}h_n(x)J(x,n,0)dxdn+\int_{0}^{\infty} \int^{\i}_{0}dt \int_{\mathbb{R}^{d}} 
f(x,t,dn)\left(\frac {\partial}{\partial t}+ d(n)
\triangle_x\right) J(x,n,t)
\\  && +\int^{\i}_{0}  \int _{L^{-1}}^{L} \int _{L^{-1}}^{L} \int_{\mathbb{R}^{d}}  \b(m,n) f(x,t,dm) f(x,t,dn)  
\tilde J(x,m,n,t)  dx dt+Err^4(L)=0.
\end{eqnarray}
It remains to replace $L^{-1}$ and $L$ with $0$ and $\i$ respectively. For this, recall that by assumption, there exists $\ell$
 such that $J(x,m,t)=0$ if $m\notin (\ell^{-1},\ell)$. Hence, when $\tilde J(x,m,n,t) \neq 0$, we must have that
$m+n> \ell^{-1}$ and  $\min\{m,n\}< \ell$. By the first remark we made after the statement of Theorem 1.1, we know that $\b\le\a$.
From the second part of Hypothesis 1.1 we deduce that there exists a constant $c_5=c_5(\ell)$ 
such that $\b(m,n)\le\a(m,n)\le c_5\g_2(m)\g_2(n)$ provided that $m+n> \ell^{-1}$ and  $\min\{m,n\}< \ell$. (Here we are using the fact
that $d(m)^{d/2}\phi^{d-1}$ is uniformly positive and bounded over the interval $[\ell^{-1}/2,\ell]$.)
 On the other hand, we know by part 1 of Theorem 5.1,
\[
\int_0^T\int_{|x|\le L_1}\int_0^\i\int_0^\i \g_2(n)\g_2(m)f(x,t,dm)f(x,t,dn)dxdt<\i,
\]
$\cP$-almost surely, where $L_1$ is chosen so that the set $\{|x|\le L_1\}$ contains the support of $\tilde J$ in the spatial variable.
 From this we deduce 
\begin{eqnarray*}
&&\lim_{L\to\i}\int_0^T\int\int_0^\i\int_0^\i \b(m,n)f(x,t,dm)f(x,t,dn)\\
&&\quad\quad \quad 1\!\!1\left(\max\{m,n\}\ge
L{\text{ or }}\min\{m,n\}\le L^{-1}\right)\tilde J(x,m,n,t)dxdt=0.
\end{eqnarray*}
This allows us to replace $L^{-1}$ and $L$ with $0$ and $\i$ respectively in (6.15), concluding that $f(x,t,dn)$ solves (1.1) 
weakly in the sense of (1.6).
\qed

As we stated in Section~1, the family $\cP^\e$ is defined on a compact metric space $\cX$ which consists of measures $\mu(dx,dn,dt)$
which are absolutely continuous with respect to the time variable. This can be proved by standard arguments. 

\begin{lemma}
Every measure $\mu\in \cX$ is of the form $\mu(dx,dn,dt)=g(dx,dn,t)dt$.
\end{lemma}

{\bf {Proof.}} Let $J_k :\bR^d\times [0,\i)\to \bR$, $k\in \bN$ be a sequence of linearly independent continuous functions of
 compact support such that $J_1=1$ and the linear span 
$Y$ of this sequence is dense in the space of continuous functions of compact support. 
Given $\mu\in \cX$, it is not hard to show that for each $k$, 
there exists a measurable function $G_{J_k}:[0,T]\to\bR$ such that $\|G_{J_k}\|_{L^\i}\le Z\sup_{x,n}|J_k(x,n)|$, and 
\[
\int_{\bR^d}\int_0^\i J_k(x,n)\mu(dx,dn,dt)=G_{J_k}(t)dt.
\]
We wish to define $G_J$ for every continuous $J$ of compact support. Note that each $G_{J_k}$ is defined almost everywhere in the interval
$[0,\i)$.
 For our purposes, we need to construct $G_J$ in such a way that for almost all $t$, the operator $J\mapsto G_J(t)$ is linear.
 For this, let us set
$G_J=r_1G_{J_1}+\dots+r_lG_{J_l}$ when $J=r_1J_1+\dots+r_lJ_l$ with $r_1,\dots,r_l$ rational. The set of such $J$ is denoted by $Y'$. 
Since $Y'$ is countable, There exists a set $A\subset [0,\i)$ of $0$ Lebesgue measure, such that for $t\notin A$, 
the operator $J\mapsto G_J(t)$ from $Y'$ to $\R$ is linear over rationals. By denseness of rationals, we can extend  $J\mapsto
G_J(t)$ for $J\in Y$ and $t\notin A$.
For such $(J,t)$, 
\[
\int_{\bR^d}\int_0^\i J(x,n)\mu(dx,dn,dt)=G_{J}(t)dt.
\]
We then take a point in $[0,\i)-A$ and use Riesz Representation Theorem to find a measure $g(dx,dn,t)$ such that 
\[
G_J(t)=\int_{\bR^d}\int_0^\i J(x,n)g(dx,dn,t),
\]
for every $J\in Y$. Hence
\[
\int_{\bR^d}\int_0^\i J(x,n)\mu(dx,dn,dt)=
\int_{\bR^d}\int_0^\i J(x,n)g(dx,dn,t)dt.
\]
for every $J\in Y$. This completes the proof.
\qed

\end{section}

\maketitle \setcounter{equation}{0}
\begin{section}{Entropy}

In this section, we establish entropy-like inequalities to show that the macroscopic density $g$ is
 absolutely continuous with respect to Lebesgue measure.

\noindent{\bf Proof of Theorem 1.2.} \\ 
\noindent{\bf Step1.} Recall that initially we have ${\cN}$ particles.
 We choose $I_{\mathbf q(0)}=\{1,\dots,{\cN}\}$, and label the initial particles as
$(x_1,m_1),\dots,(x_{\cN},m_{\cN})$. If a coagulation occurs at time $t$, 
one of 
the coagulating particles disappears from the system, and $I_{\bf q} \subseteq \big\{ 1,\ldots, \cN \big\}$
 satisfies $\big\vert I_{{\bf q}(t^+)}\big\vert =\big\vert I_{{\bf q}(t)}\big\vert - 1$. 
We write ${\cN}(\mathbf q)=|I_{\mathbf q}|$ for the  number of particles of the configuration $\mathbf q$.
Note that  ${\cN}(\mathbf q)$ takes values in the set
$\{1,\dots,{\cN}\}$. We write $F(\mathbf q,t)\nu_N(d\mathbf q)$ for the law of 
$\mathbf q(t),$ and define
\[
H_N(t)=\int F(\mathbf q,t)\log F(\mathbf q,t)\ \nu_N(d\mathbf q).
\]
By standard arguments,  
\begin{equation}
\frac{\partial H_N}{\partial t}(t)= \int \Big( \bL  (\log F)(\mathbf q,t) \Big) F(\mathbf q,t)\nu_N(d\mathbf q)=\O_1+\O_2,
\end{equation}
where 
\begin{eqnarray*}
\O_1&=&\int \Big( \mathbb A_0(\log F)(\mathbf q,t) \Big) F(\mathbf q,t)\nu_N(d\mathbf q), \\
\O_2&=& \int \Big( \mathbb A_c(\log F)(\mathbf q,t) \Big)  F(\mathbf q,t)\nu_N(d\mathbf q).
 \end{eqnarray*}
We have
\begin{eqnarray*}
\O_1&=&\int\sum_{i\in I_{\mathbf q}}d(m_i) \big( \triangle_{x_i}F \big) \log F \ d\nu_N\\
&=&-\int\sum_{i\in I_{\mathbf q}}d(m_i)\frac{|\nabla_{x_i}F|^2}{F}d\nu_N+\int\sum_{i\in I_{\mathbf q}} d(m_i) \nabla_{x_i}F\cdot x_i\ d\nu_N\\
&=&-\int\sum_{i\in I_{\mathbf q}}d(m_i)\frac{|\nabla_{x_i}F|^2}{F}d\nu_N- \int\sum_{i\in I_{\mathbf q}}d(m_i)(d-|x_i|^2)F\ d\nu_N\\
&\le  & D\int\sum_{i\in I_{\mathbf q}}|x_i|^2F\ d\nu_N,
\end{eqnarray*}
where we integrated by parts for the second and third equality, and $D$ is an upper bound for the function $d(\cdot)$. 
To bound the right-hand side, we use the Markov property of the process $\mathbf q(t)$ to write
\begin{eqnarray}
\bE_N\sum_{i\in I_{\mathbf q(t)}}|x_i(t)|^2 & \leq & \bE_N\sum_{i\in I_{\mathbf q(0)}}|x_i(0)|^2+2d\int_0^t\bE_N\sum_{i\in I_{\mathbf q(s)}}
d(m_i(s))ds \nonumber \\
 & \leq & c \epsilon^{2 - d} + 2 d tD Z \epsilon^{2-d}, \nonumber 
\end{eqnarray}
where, in the first inequality, we used that the coagulation is non-positive,
 which follows from our assumption that a particle, newly born in a coagulation event, 
is placed in the location of one of the departing particles. The second inequality is due to our assumption that 
$D$ is a uniform upper bound on $d:(0,\infty) \to (0,\infty)$ and to the hypothesis we make on the initial condition.      
We learn that
\begin{equation}
\O_1\le c_1 (t+1) \e^{2-d}.
\end{equation}

We now concentrate on the contribution coming from 
coagulations, namely the expression $\O_2$.
This expression equals
\begin{eqnarray*}
&& \int \sum_{i,j\in I_{{\mathbf q}}}V_\e(x_i-x_j) \a(m_i,m_j)\left[\frac{m_i}{m_i+m_j}\log\frac{F(S^1_{i,j}{\mathbf q},t)}{F({\mathbf q},t)}+
\frac{m_j}{m_i+m_j}\log\frac{F(S^2_{i,j}{\mathbf q},t)}{F({\mathbf q},t)}\right]F({\mathbf q},t)\ \nu_N(d{\mathbf q})\\
&&\le \int \sum_{i,j\in I_{{\mathbf q}}}V_\e(x_i-x_j) \a(m_i,m_j)\left[\frac{m_i}{m_i+m_j}{F(S^1_{i,j}{\mathbf q},t)}+
\frac{m_j}{m_i+m_j}{F(S^2_{i,j}{\mathbf q},t)}\right]\ \nu_N(d{\mathbf q})\\
&&= \int \sum_{i,j\in I_{{\mathbf q}}}V_\e(x_i-x_j) \a(m_i,m_j)F(S^1_{i,j}{\mathbf q},t)\ \nu_N(d{\mathbf q}), 
\end{eqnarray*}
where we used the elementary inequality $\log x\le x$ for the second line. To bound this, we first observe 
\[
\int V_\e(x_i-x_j)(2\pi)^{-d/2}\exp\left(-\frac{|x_j|^2}2\right)dx_i\le (2\pi)^{-d/2} \int V_\e(x_i-x_j)dx_i\le C \e^{d-2}.
\]
We then make a change of variables $m_i+m_j\mapsto m_i$. As a result, $\O_2$ is bounded above by
\[
\e^{d-2} \int \sum_{i\in I_{ \mathbf q}}\rho(m_i)F(\mathbf q,t)d\nu_N(d\mathbf q),
\]
where the function $\rho$ is defined (\ref{defrho}).

From the second part of Hypothesis 1.3, we deduce that
 $\O_2$ is bounded by a constant multiple of $\e^{d-2}$. This, the first part of Hypothesis 1.3, and (7.2) yield
\begin{equation}
H_N(t)\le c_2 (t+1) \e^{d-2}.
\end{equation}

\noindent{\bf Step 2.} Note that by Sanov's theorem, the empirical measure $\e^{d-2}\sum_i\d_{(x_i,m_i)}$ satisfies
 a large deviation principle with 
respect to the measure $\nu_N$ as $\e\to 0$.
The large deviation rate function  $\cI(g)=\i$ unless $g(dx,dn)=f(x,n)r(x,n)dxdn$ and if such a function $f$ exists, then 
\[
\cI(g)=\int_0^\i\int (f\log {f}-f+1){r}\ dxdn.
\]
By an argument similar to the proof of Lemma 6.3 of \cite{GPV}, 
we can use (7.3) to deduce that if $\cP$ is any limit point of the sequence $\cP^\e$, then 
\[
\int \cI(g(\cdot,t))\ \cP(dg)<\i,
\]
for every $t$. This completes the proof of Theorem 1.2.
\qed

\end{section}
\end{section}
\begin{section}{Appendix: Scaling of the continuous Smoluchowski equation}
We comment on the scaling satisfied by the system (\ref{eq1.1}), under the assumptions that
$$
 d(n) = n^{-\phi}
$$
and
\begin{equation}\label{betabd}
\beta(n,m) = n^{\eta} + m^{\eta},
\end{equation}
with $\phi,\eta \in [0,\infty)$.
Rescaling the equations,
\begin{equation}\label{gneqn}
g_n(x,t) = \lambda^{\alpha} f_{n \lambda^{\gamma}}\big( \lambda^{\tau}x,\lambda t \big),
\end{equation}
we note that $g_n$ satisfies (\ref{eq1.1}) provided that
\begin{equation}\label{condfree}
 1 - \gamma \phi - 2 \tau = 0 
\end{equation}
and
\begin{equation}\label{condint}
 - \alpha  + \gamma \big( 1 + \eta  \big) + 1 = 0, 
\end{equation}
(\ref{condfree}) ensuring that the free motion term is preserved, (\ref{condint}) the interaction term.
The mass 
$$
h_f(t) = \int_0^{\infty} n \int_{\mathbb{R}^d} f_n \big( x,t \big) dx\ dn,
$$
which, formally at least, is conserved in time, is mapped by the rescaling to 
\begin{equation}\label{hgeqn}
 h_g(t) = \lambda^{\alpha - \tau d - 2 \gamma} h_f(\l t). 
\end{equation}
The mass, then, is conserved by the rescaling 
provided that 
\begin{equation}\label{condener}
  \alpha - \tau d - 2 \gamma = 0.
\end{equation}
In the critical case, where each of (\ref{condfree}),
(\ref{condint}) and (\ref{condener}) is satisfied,  
 we have that
\begin{displaymath}
   \gamma = \frac{d/2 - 1}{\eta + \phi d/2 - 1},
\end{displaymath}
\begin{displaymath}
 \alpha = \frac{d/2 \big( \phi + \eta + 1 \big) - 2}{\eta + \phi d/2 - 1}
\end{displaymath}
and
\begin{equation}\label{eqnchi}
\tau = \frac{\eta + \phi - 1}{2 \big( \eta + \phi d/2 - 1 \big)}.
\end{equation}
In the case that the dimension $d=2$, the values $\gamma = 0$, $\alpha = 1$ and $\tau = 1/2$ are adopted, whatever the values taken for 
the input parameters $\phi$ and $\eta$. The only critical scaling, then, leaves the mass unchanged and performs a diffusive rescaling of 
space-time.
 
Regarding the critical scaling, we recall from Remark 1.2 of \cite{HR3} that 
the condition $\eta + \phi = 1$, which is a natural transition 
for the rescaling $g_n$ (as is apparent from (\ref{eqnchi})), represents the limit of the parameter range for which uniqueness and 
mass-conservation of the solution of (\ref{eq1.1}) are proved: indeed, the condition required by \cite{HR3} is $\eta + \phi < 1$, 
along with some hypothesis on the initial data.  

Do we expect the complementary condition $\eta + \phi \geq 1$ to have physical meaning? To consider this question, we take positive and 
fixed $\phi$ and $\eta$, and consider the rescaling (\ref{gneqn}) under the constraints (\ref{condfree}) and (\ref{condint}). Seeking to 
understand the formation of massive particles, rather than spatial blow-up, we fix $\tau = 0$. We are led to
\begin{equation}\label{gpeqn}
 \gamma = \phi^{-1}
\end{equation}
and 
\begin{equation}\label{alpheqn}
\alpha = 1 + \frac{1 + \eta}{\phi}.
\end{equation}
Returning to (\ref{gneqn}), a self-similar blow-up profile is consistent with the scaling 
$$ 
t^{-\alpha} f_{n t^{-\gamma}} \big( x , 1 \big)  
$$
given by $\lambda = t^{-1}$ provided that its mass (\ref{hgeqn}) does not grow to infinity as $\lambda \to 0$. We have set $\tau = 0$: 
as such, the condition that ensures this is $\alpha - 2 \gamma \geq 0$, which, by (\ref{gpeqn}) and (\ref{alpheqn}), amounts to the 
inequality $\phi + \eta \geq 1$.

We conclude that considerations of scaling would in principle permit a blow-up in the equations in the mass variable under 
the condition that $\eta + \phi \geq 1$. The blow-up we considered is in a low $\lambda$ limit, which corresponds to heavy mass
 at late times: as such, it should be considered not as a gelation, in which particles of infinite mass develop in finite time,
 but rather as the appearance of populations of arbitrarily heavy particles at correspondingly high time-scales. Expressed more
 precisely, the weak form of blow-up considered is the statement that, for each $K$ strictly less than the total initial mass 
$\int_0^{\infty} \int_{\mathbb{R}^d} mf_m(x,0) dx dm$ and any $m_0 \in \mathbb{R}^+$, there exists $t \in [0,\infty)$, 
\begin{equation}\label{massfin}
 \int_{m_0}^{\infty} \int_{\mathbb{R}^d} m f_m(x,t) dx dm > K.
\end{equation}
(This condition is correct in the absence of gelation. Gelation would remove mass from all finite levels. Note also that the
 absence of fragmentation in (\ref{eq1.1}) means that, in fact, (\ref{massfin}) implies the stronger statement that most of 
the mass accumulates in arbitrarily high levels at all sufficiently late times.) In dimension $d \geq 3$, (1.11)
 of Theorem 1.1 in \cite{HR3} shows that the discrete analogue of (\ref{massfin}) fails if $\eta + \phi < 1$. 

A parallel may be drawn between the Smoluchowski PDE and the non-linear Schr\"odinger equation.
Consider, for example, a solution of cubic defocussing NLS,  
$u:\mathbb{R}^d \times \mathbb{R}^+ \to \mathbb{C}$ of 
\begin{equation}\label{nls}
i \frac{\partial}{\partial t} u  - \Delta u  = - \vert u \vert^2 u,
\end{equation}
may be written in Fourier space as
\begin{equation}\label{nlsfour}
i \frac{\partial}{\partial t} \hat{u}  - \vert \xi \vert^2 \hat{u}  = 
- \int \int \hat{u} (\xi - \eta) \hat{u}(\sigma) \hat{\overline{u}}(\eta - \sigma)  d \eta d \sigma.
\end{equation}
We see that the mass variable in (\ref{eq1.1}) may be viewed as analogous to the frequency variable in (\ref{nlsfour}): 
the non-linear interaction term in each case is a type of convolution. Pursuing the analogy,
 the quantity $\frac{1}{2} \vert\vert \nabla u \vert\vert_{2}^2 + \frac{1}{4} \vert\vert  u \vert\vert_{4}^4$ is formally conserved in NLS, 
as is the mass $\int_0^{\infty} \int_{\mathbb{R}^d} m f_m dx dm$ for the Smoluchowski PDE. For NLS, the term weak turbulence refers to the
 growth to infinity in time of the $H^s$ norm
$$
 \vert \vert u \vert \vert_{H^s} = \int \vert \hat{u}(\xi) \vert^2 \vert \xi \vert^{2s} d \xi,
$$
for some $s > 1$, a circumstance that is anticipated in (\ref{nls}) 
in a periodic domain. (See Section $II.2$ of \cite{Bourgain} 
for a discussion.)  
The counterpart of weak turbulence for the system (\ref{eq1.1}) is
$$
\int_0^{\infty} \int_{\mathbb{R}^d}{m^{r} f_m(x,t)} dx dm \to \infty \, \, \textrm{as $t \to \infty$,}
$$
for some $r > 1$. (Note that (\ref{massfin}) implies this statement for every $r > 1$ on a subsequence of times.)

Comparing the system (\ref{eq1.1}) to its spatially homogeneous counterpart, given in the discrete case by 
$\big\{ a_n:[0,\infty) \to [0,\infty) : n \in \mathbb{N} \big\}$ satisfying
\begin{equation}\label{homsys}
 \frac{d}{dt}a_n(t) = \sum_{m=1}^{n-1} \beta(m,n-m) a_m(t) a_{n-m}(t)  - 2 \sum_{m=1}^{\infty} \beta (m,n) a_m(t) a_n(t),
\end{equation}
we see the stabilizing role of diffusion: for example, it is easy to see that, taking $\beta(n,m)$ identically equal to a constant in 
(\ref{homsys}) ensures 
the analogue of (\ref{massfin}), 
while we have seen in the spatial case 
that scaling arguments do not disallow (\ref{massfin}) under the condition that   $\eta + \phi \geq 1$.

Regarding the prospect of proving mass-conservation for at least some part of the parameter space where $\phi + \eta \geq 1$, 
we comment that, in \cite{HR3}, hypotheses of the form $\beta(n,m) \leq n^{\eta} + m^{\eta}$ were used. It may be that, 
if $\beta(n,m) \leq n^{1 + \epsilon} + m^{1 + \epsilon}$ or $\beta(n,m) \leq n^{1/2 + \epsilon} m^{1/2 + \epsilon}$ (with $\epsilon > 0$ 
a small constant), but $\beta$ is permitted to have space-time dependence subject to such a bound, then gelation is more liable to occur.
 As such, an argument for mass-conservation  would have to exploit the assumption that $\beta(n,m)$ is constant in space-time,
 in a way that those in \cite{HR3} did not.   

\end{section}
\bibliography{biblion}

\end{document}